\documentclass[11pt]{article}
\usepackage{amssymb}
\usepackage{amsmath}
\usepackage{times}
\usepackage{graphicx}

\title{Baire-class $\xi$ colorings: the first three levels\indent}
\author{Dominique LECOMTE and Miroslav ZELENY$^{1,2}$}
\date{\today}

\def\ufootnote#1{\let\savedthfn\thefootnote\let\thefootnote\relax
\footnote{#1}\let\thefootnote\savedthfn\addtocounter{footnote}{-1}}

\newcommand{\Ana}{{\it\Sigma}^{1}_{1}}

\newcommand{\Ca}{{\it\Pi}^{1}_{1}}
\newcommand{\Boraone}{{\it\Sigma}^{0}_{1}}

\newcommand{\Borel}{{\it\Delta}^{1}_{1}}

\newcommand{\borel}{{\bf\Delta}^{1}_{1}}

\newcommand{\boraone}{{\bf\Sigma}^{0}_{1}}
\newcommand{\boratwo}{{\bf\Sigma}^{0}_{2}}
\newcommand{\borathree}{{\bf\Sigma}^{0}_{3}}

\newcommand{\boraxi}{{\bf\Sigma}^{0}_{\xi}}

\newcommand{\borone}{{\bf\Delta}^{0}_{1}}
\newcommand{\bortwo}{{\bf\Delta}^{0}_{2}}
\newcommand{\borthree}{{\bf\Delta}^{0}_{3}}
\newcommand{\borfour}{{\bf\Delta}^{0}_{4}}
 
\newcommand{\bormone}{{\bf\Pi}^{0}_{1}}

\newcommand{\bormtwo}{{\bf\Pi}^{0}_{2}}
\newcommand{\bormthree}{{\bf\Pi}^{0}_{3}}

\newcommand{\bormlxi}{{\bf\Pi}^{0}_{<\xi}}

\newcommand{\bormxi}{{\bf\Pi}^{0}_{\xi}}

\newcommand{\borxi}{{\bf\Delta}^{0}_{\xi}}
\newcommand{\borme}{{\bf\Pi}^{0}_{\eta}}

\newtheorem{thm} {Theorem} [section]
\newtheorem{defi} [thm] {Definition}
\newtheorem{cor} [thm] {Corollary}
\newtheorem{lem} [thm] {Lemma}
\newtheorem{prop} [thm] {Proposition}

\begin{document}

\maketitle

\centerline{$\bullet$ Universit\' e Paris 6, Institut de Math\'ematiques de Jussieu, Projet Analyse Fonctionnelle}

\centerline{Couloir 16-26, 4\`eme \'etage, Case 247, 4, place Jussieu, 75 252 Paris Cedex 05, France}

\centerline{dominique.lecomte@upmc.fr}\bigskip

\centerline{$\bullet$ Universit\'e de Picardie, I.U.T. de l'Oise, site de Creil,}

\centerline{13, all\'ee de la fa\"\i encerie, 60 107 Creil, France}\bigskip

\centerline{$\bullet^1$ Charles University, Faculty of Mathematics and Physics, Department of Mathematical Analysis}

\centerline{Sokolovsk\'a 83, 186 75 Prague, Czech Republic}

\centerline{zeleny@karlin.mff.cuni.cz}\bigskip\bigskip\bigskip\bigskip\bigskip\bigskip

\ufootnote{{\it 2010 Mathematics Subject Classification.}~Primary: 03E15, Secondary: 54H05}

\ufootnote{{\it Keywords and phrases.}~Borel chromatic number, Borel class, coloring,  dichotomy, Hurewicz, partition, product}

\ufootnote{$^2$ The work is a part of the research project MSM 0021620839 financed by MSMT and partly supported by the grant GA\v CR 201/09/0067.}

\ufootnote{\bf Acknowledgements.\rm\ This work started in summer 2009, and was continued during  spring 2010, when the first author was invited at the University of Prague by the second one. The first author is very grateful to the second one for that. We would also like to thank Alain Louveau for some nice remarks made in the Descriptive Set Theory Seminar of the University Paris 6. They were useful to simplify our presentation.}

\noindent {\bf Abstract.} The $\mathbb{G}_0$-dichotomy due to Kechris, Solecki and Todor\v cevi\'c characterizes the analytic relations having a Borel-measurable countable coloring. We give a version of the $\mathbb{G}_0$-dichotomy for $\boraxi$-measurable countable colorings when $\xi\!\leq\! 3$. A 
$\boraxi$-measurable countable coloring gives a covering of the diagonal consisting of countably many 
$\boraxi$ squares. This leads to the study of countable unions of $\boraxi$ rectangles. We also give a Hurewicz-like dichotomy for such countable unions when $\xi\!\leq\! 2$.

\vfill\eject

\section{$\!\!\!\!\!\!$ Introduction}\indent

 The reader should see [K] for the standard descriptive set theoretic notation used in this paper. We study a definable coloring problem. We will need some more notation:\bigskip
 
\noindent\bf Notation.\rm ~The letters $X$, $Y$ will refer to some sets. We set 
$\Delta (X)\! :=\!\{ (x_0,x_1)\!\in\! X^2\mid x_0\! =\! x_1\}$.

\begin{defi} (1) Let $A\!\subseteq\! X^2$. We say that $A$ is a $digraph$ if $A$ does not meet $\Delta (X)$.
\smallskip

\noindent (2) Let $A$ be a digraph. A $countable~coloring$ of $(X,A)$ is a map 
$c\! :\! X\!\rightarrow\!\omega$ such that $A$ does not meet $(c\!\times\! c)^{-1}\big(\Delta (\omega )\big)$.\end{defi}

 In [K-S-T], the authors characterize the analytic digraphs of having a Borel countable coloring. The characterization is given in terms of the following notion of comparison between relations.\bigskip

\noindent\bf Notation.\rm ~Let $X,Y$ be Polish spaces, $A$ (resp., $B$) be a relation on $X$ (resp., $Y$), and $\bf\Gamma$ be a class of sets. We set
$$(X,A)\preceq_{\bf\Gamma}(Y,B)~\Leftrightarrow ~\exists f\! :\! X\!\rightarrow\! Y~~{\bf\Gamma}\mbox{-measurable with }A\!\subseteq\! (f\!\times\! f)^{-1}(B).$$
In this case, we say that $f$ is a ${\bf\Gamma}$-$measurable\ homomorphism$ from $(X,A)$ into $(Y,B)$. This notion essentially makes sense for digraphs (we can take $f$ to be constant if $B$ is not a digraph).\bigskip

 We also have to introduce a minimum digraph without Borel countable coloring:\bigskip

\noindent $\bullet$ Let $\psi\! :\!\omega\!\rightarrow\! 2^{<\omega}$ be a natural bijection. More specifically, $\psi (0)\! :=\!\emptyset$ is the sequence of length $0$, $\psi (1)\! :=\! 0$, $\psi (2)\! :=\! 1$ are the sequences of length $1$, and so on. Note that $|\psi (n)|\!\leq\! n$ if $n\!\in\!\omega$. Let 
$n\!\in\!\omega$. As $|\psi (n)|\!\leq\! n$, we can define $s_n\! :=\!\psi (n)0^{n-|\psi (n)|}$. The crucial properties of the sequence $(s_{n})_{n\in\omega}$ are the following:\bigskip

- $(s_n)_{n\in\omega}$ is $dense$ in $2^{<\omega}$. This means that for each $s\!\in\! 2^{<\omega}$, there is $n\!\in\!\omega$ such that $s_n$ extends $s$ (denoted $s\!\subseteq\! s_n$).\smallskip

- $|s_n|\! =\! n$.\bigskip

\noindent $\bullet$ We put
$\mathbb{G}_0\! :=\!\{ (s_n0\gamma ,s_n1\gamma )\mid n\!\in\!{\omega}\mbox{ and }
\gamma\!\in\! 2^{\omega}\}\!\subseteq\! 2^\omega\!\times\! 2^\omega$. Note that $\mathbb{G}_0$ is analytic (in fact difference of two closed sets) since the map $(n,\gamma )\!\mapsto\! (s_n0\gamma ,s_n1\gamma )$ is continuous.\bigskip

 The previous definitions were given, when ${\bf\Gamma}\! =\!\borel$, in [K-S-T], where the following is proved:

\begin{thm} (Kechris, Solecki, Todor\v cevi\'c) Let $X$ be a Polish space, and $A$ be an analytic relation on $X$. Then exactly one of the following holds:\smallskip  

\noindent (a) There is a Borel countable coloring of $(X,A)$, i.e., 
$(X,A)\preceq_{\borel}\big(\omega ,\neg\Delta (\omega )\big)$,\smallskip  

\noindent (b) $(2^{\omega},\mathbb{G}_0)\preceq_{\boraone}(X,A)$.\end{thm} 

 This result had several developments during the last decade. Here is a non-exhaustive list:\bigskip
 
\noindent - We can characterize the potentially open sets via a Hurewicz-like test, and in finite dimension it is a consequence of the previous result. Let us specify this. The following definition can be found in [Lo2] (see Definition 3.3).
 
\begin{defi} (Louveau) Let $X,Y$ be Polish spaces, $A$ be a Borel subset of $X\!\times\! Y$, and 
$\bf\Gamma$ be a Borel class. We say that $A$ is $potentially~in~{\bf\Gamma}$ $\big($denoted 
$A\!\in\!\mbox{pot}({\bf\Gamma})\big)$ if we can find a finer Polish topology $\sigma$ $($resp., $\tau )$ on $X$ $($resp., $Y)$ such that $A\!\in\! {\bf\Gamma}\big( (X, \sigma )\!\times\! (Y, \tau )\big)$.\end{defi}

 The $\mbox{pot}(\boraone )$ sets are the countable unions of Borel rectangles. A consequence of this is that the Borel hierarchy built on the Borel rectangles is exactly the hierarchy of the classes of the sets potentially in some Borel class. The good notion of comparison to study the 
$\mbox{pot}({\bf\Gamma})$ sets is as follows (see [L3]). Let $X_0,X_1,Y_0,Y_1$ be Polish spaces, and $A^\varepsilon_0,A^\varepsilon_1\!\subseteq\! X_\varepsilon\!\times\! Y_\varepsilon$ be disjoint. We set\bigskip

\leftline{$(X_0,Y_0,A^0_0,A^0_1)\leq (X_1,Y_1,A^1_0,A^1_1)\Leftrightarrow$}\smallskip

\rightline{$\exists f\! :\! X_0\!\rightarrow\! X_1~~\exists g\! :\! Y_0\!\rightarrow\! Y_1\mbox{ continuous with }
A^0_\varepsilon\!\subseteq\! (f\!\times\! g)^{-1}(A^1_\varepsilon )\mbox{ for each }\varepsilon\!\in\! 2.$}\bigskip

 The following theorem is proved in [L1], and is a consequence of Theorem 1.2:
 
\begin{thm} Let $X,Y$ be Polish spaces, and $A_0,A_1$ be disjoint analytic subsets of $X\!\times\! Y$. Then exactly one of the following holds:\smallskip

\noindent (a) The set $A_0$ can be separated from $A_1$ by a 
$\mbox{pot}(\boraone )\! =\! (\borel\!\times\!\borel )_\sigma$ set (i.e., there is 
$S\!\in\!\mbox{pot}(\boraone )$ with $A_0\!\subseteq\! S\!\subseteq\!\neg A_1$),\smallskip

\noindent (b) $\big( 2^\omega ,2^\omega ,\Delta (2^\omega ),\mathbb{G}_0\big)\leq (X,Y,A_0,A_1)$.
\end{thm}
 
 In [L1], it is also proved that we cannot have $f$ one-to-one in Theorem 1.2.(b) in general. It is easy to check that Theorem 1.2 is also an easy consequence of Theorem 1.4. This means that the study of the Borel countable colorings is highly related to the study of countable unions of Borel rectangles.\bigskip
 
\noindent - We can extend Theorem 1.2 to any finite dimension, and also in infinite dimension if we change the space in which lives the infinite dimensional version of $\mathbb{G}_0$ (see [L2]).\bigskip
 
\noindent - B. Miller recently developped some techniques to recover many dichotomy results of descriptive set theory, but without using effective descriptive set theory. He replaces it with some versions of Theorem 1.2. In particular, he can prove Theorem 1.2 without effective descriptive set theory.\bigskip

  When $A$ is Borel, it is natural to ask about the relation between the Borel class of $A$ and that of the coloring $f$ when Theorem 1.2.(a) holds. This leads to consider $\borxi$-measurable countable colorings (or equivalently $\boraxi$-measurable countable colorings). We have the following conjecture:\bigskip
   
\noindent\bf Conjecture 1\it\ Let $1\!\leq\!\xi\! <\!\omega_1$. Then there are\smallskip

\noindent - a $0$-dimensional Polish space $\mathbb{X}_\xi$,\smallskip

\noindent - an analytic relation $\mathbb{A}_{\xi}$ on $\mathbb{X}_\xi$\smallskip

\noindent such that for any $0$-dimensional Polish space $X$, and for any analytic relation $A$ on $X$, exactly one of the following holds:\smallskip  

\noindent (a) $(X,A)\preceq_{\borxi}\big(\omega ,\neg\Delta (\omega )\big)$,\smallskip  

\noindent (b) $(\mathbb{X}_\xi ,\mathbb{A}_{\xi})\preceq_{\boraone}(X,A)$.\rm\bigskip

 We will prove it when $1\!\leq\!\xi\!\leq\! 2$, and in these cases we do not have to assume that $A$ is analytic. We will also prove it when $\xi\! =\! 3$, which is much more difficult. We should not have to assume that $X$ is $0$-dimensional when $\xi\!\geq\! 2$, but we have to do it when $\xi\! =\! 1$.
 
\vfill\eject

 We saw that the study of the Borel countable colorings is highly related to the study of countable unions of Borel rectangles, and gave some motivation for studying $\boraxi$-measurable countable colorings. This motivates the study of countable unions of $\boraxi$ rectangles. Another motivation is that 
$(X,A)\preceq_{\borxi}\big(\omega ,\neg\Delta (\omega )\big)$ is equivalent to the fact that $\Delta (X)$ can be separated from $A$ by a $(\boraxi\!\times\!\boraxi )_\sigma$ set, by the generalized reduction property for the class $\boraxi$ (see 22.16 in [K]).\bigskip
   
\noindent\bf Conjecture 2\it\ Let $1\!\leq\!\xi\! <\!\omega_1$. Then there are $0$-dimensional Polish spaces $\mathbb{X}^0_\xi ,\mathbb{X}^1_\xi$ and disjoint analytic subsets 
$\mathbb{A}^0_{\xi},\mathbb{A}^1_{\xi}$ of $\mathbb{X}^0_\xi\!\times\!\mathbb{X}^1_\xi$ such that for any Polish spaces $X,Y$, and for any pair $A_0,A_1$ of disjoint analytic subsets of $X\!\times\! Y$, exactly one of the following holds:\smallskip  

\noindent (a) The set $A_0$ can be separated from $A_1$ by a $(\boraxi\!\times\!\boraxi )_\sigma$ set,\smallskip  

\noindent (b) $(\mathbb{X}^0_\xi ,\mathbb{X}^1_\xi ,\mathbb{A}^0_\xi ,\mathbb{A}^1_\xi )\leq 
(X,Y,A_0,A_1)$.\rm\bigskip

 It is trivial to prove this when $\xi\! =\! 1$. We will prove that Conjecture 2 holds when 
$\xi\!\leq\! 2$, which is significantly more and more difficult when $\xi$ increases. We use effective descriptive set theory, and give effective strengthenings of our results. The reader should see [M] for basic notions of effective descriptive set theory. In particular, we will see that to test whether an analytic relation has a $\boraxi$-measurable countable coloring, it is enough to test countably many partitions instead of continuum many. We will use the topology $T_\xi$ generated by the $\Ana\cap\bormlxi$ subsets of a recursively presented Polish space (introduced in [Lo1]) when $\xi$ is $2$ or $3$ ($T_1$ is just the basic topology). The last result can be strengthened as follows (see [L3]).

\begin{thm} Let $1\!\leq\!\xi\!\leq\! 2$. Then there are $0$-dimensional Polish spaces 
$\mathbb{X}^0_\xi ,\mathbb{X}^1_\xi$ and disjoint analytic subsets 
$\mathbb{A}^0_{\xi},\mathbb{A}^1_{\xi}$ of $\mathbb{X}^0_\xi\!\times\!\mathbb{X}^1_\xi$ such that for any recursively presented Polish spaces $X,Y$, and for any pair $A_0,A_1$ of disjoint $\Ana$ subsets of 
$X\!\times\! Y$, the following are equivalent:\smallskip

\noindent (a) The set $A_0$ cannot be separated from $A_1$ by a 
$(\boraxi\!\times\!\boraxi )_\sigma$ set.\smallskip

\noindent (b) The set $A_0$ cannot be separated from $A_1$ by a 
$\Borel\cap (\boraxi\!\times\!\boraxi )_\sigma$ set.\smallskip

\noindent (c) The set $A_0$ cannot be separated from $A_1$ by a 
$\boraone (T_\xi\!\times\! T_\xi )$ set.\smallskip

\noindent (d) $A_0\cap\overline{A_1}^{T_\xi\times T_\xi}\!\not=\!\emptyset$.\smallskip

\noindent (e) $(\mathbb{X}^0_\xi ,\mathbb{X}^1_\xi ,\mathbb{A}^0_\xi ,\mathbb{A}^1_\xi )\leq 
(X,Y,A_0,A_1)$.\end{thm}

\section{$\!\!\!\!\!\!$ Some general effective facts}\indent

 One can hope for an effective strengthening of Conjecture 1:\bigskip
    
\noindent\bf Effective conjecture 1\it\ Let $1\!\leq\!\xi\! <\!\omega_1$. We can find a $0$-dimensional Polish space $\mathbb{X}_\xi$ and an analytic relation $\mathbb{A}_{\xi}$ on $\mathbb{X}_\xi$ such that $(\mathbb{X}_\xi ,\mathbb{A}_{\xi})\not\preceq_{\borxi}\big(\omega ,\neg\Delta (\omega )\big)$, and 
for any $\alpha\!\in\!\omega^\omega$ with $1\!\leq\!\xi\! <\!\omega_1^\alpha$, for any 
$0$-dimensional recursively in $\alpha$ presented Polish space $X$, and for any $\Ana (\alpha )$ relation $A$ on $X$, one of the following holds:\smallskip  

\noindent (a) $(X,A)\preceq_{\Borel (\alpha )\cap\borxi}
\big(\omega ,\neg\Delta (\omega )\big)$,\smallskip  

\noindent (b) $(\mathbb{X}_\xi ,\mathbb{A}_{\xi})\preceq_{\boraone}(X,A)$.\rm

\vfill\eject

 We will see that this effective conjecture is true when $1\!\leq\!\xi\!\leq\! 3$. The following statement is a corollary of this effective conjecture, and is in fact a theorem: 
   
\begin{thm} Let $1\!\leq\!\xi\! <\!\omega_1^{\mbox{CK}}$, $X$ be a $0$-dimensional recursively presented Polish space, and $A$ be a $\Ana$ relation on $X$. We assume that 
$(X,A)\preceq_{\borxi}\big(\omega ,\neg\Delta (\omega )\big)$. Then 
$(X,A)\preceq_{\Borel\cap\borxi}\big(\omega ,\neg\Delta (\omega )\big)$.\end{thm}

 A consequence of this is that to test whether an analytic relation has a $\boraxi$-measurable countable coloring, it is enough to test countably many partitions instead of continuum many. Another consequence is the equivalence between Conjecture 1 and the Effective conjecture 1. We have in fact preliminary results that will help us to prove also the equivalence between (a)-(d) in Theorem 1.5, in the general case. 
 
\begin{lem} Let $1\!\leq\!\xi\! <\!\omega_1^{\mbox{CK}}$, $X, Y$ be recursively presented Polish spaces, 
$A\!\in\!\Ana (X)\cap\boraxi$, $B\!\in\!\Ana (Y)\cap\boraxi$, and $C\!\in\!\Ana (X\!\times\! Y)$ disjoint from 
$A\!\times\! B$. Then there are $A',B'\!\in\!\Borel\cap\boraxi$ such that $A'\!\times\! B'$ separates 
$A\!\times\! B$ from $C$. This also holds for $\bormxi$ instead of $\boraxi$.\end{lem}
 
\noindent\bf Proof.\rm\ Note that $A$ and $\{ x\!\in\! X\mid\exists y\!\in\! B~~(x,y)\!\in\! C\}$ are disjoint 
$\Ana$ sets, separable by a $\boraxi$ subset of $X$. By Theorems 1.A and 1.B in [Lo1], there is 
$A'\!\in\!\Borel\cap\boraxi$ separating these two sets. Similarly, $B$ and 
$\{ y\!\in\! Y\mid\exists x\!\in\! A'~~(x,y)\!\in\! C\}$ are disjoint $\Ana$ sets, and there is 
$B'\!\in\!\Borel\cap\boraxi$ separating these two sets. The proof for $\bormxi$ is identical to the one for 
$\boraxi$.\hfill{$\square$}

\begin{thm} Let $1\!\leq\!\xi\! <\!\omega_1^{\mbox{CK}}$, $X, Y$ be recursively presented Polish spaces, and $A_0,A_1$ be disjoint $\Ana$ subsets of $X\!\times\! Y$. We assume that $A_0$ is separable from $A_1$ by a $\big(\boraxi\!\times\!\boraxi\big)_\sigma$ set. Then $A_0$ is separable from $A_1$ by a 
$\Borel\cap\big( (\Borel\cap\boraxi )\!\times\! (\Borel\cap\boraxi)\big)_\sigma$ set.\end{thm}

\noindent\bf Proof.\rm\ By Example 2 of Chapter 3 in [Lo2], the family $\big( N(n,X)\big)_{n\in\omega}$ is regular without parameter. By Corollary 2.10 in [Lo2], $\bormxi (X)$, as well as 
$\boraxi (X)\! =\!\big(\bigcup_{\eta <\xi}~\borme (X)\big)_\sigma$,  are regular without parameter. By Theorem 2.12 in [Lo2], $\boraxi (X)\!\times\!\boraxi (Y)$ is also regular without parameter. By Theorem 2.8 in [Lo2], the family $\Phi\! :=\!\big(\boraxi (X)\!\times\!\boraxi (Y)\big)_\sigma$ is separating, which implies the existence of $S\!\in\!\Borel\cap\Phi$ separating $A_0$ from $A_1$.\bigskip

 With the notation of [Lo2], let $n$ be an integer with $(0^\infty ,n)\!\in\! W$ and $C_{0^\infty ,n}\! =\! S$. Then $(0^\infty ,n)$ is in $W_\Phi$, which by Theorem 2.8.(ii) in [Lo2] is 
$$\left\{ (\alpha ,n)\!\in\! W\mid\exists\beta\!\in\!\Borel (\alpha )~~\forall m\!\in\!\omega ~~
\big(\alpha ,\beta (m)\big)\!\in\! W_{\boraxi (X)\times\boraxi (Y)}\mbox{ and }
C_{\alpha ,n}\! =\!\bigcup_{m\in\omega}~C_{\alpha ,\beta (m)}\right\}.$$
This implies that $S\!\in\!\Borel\cap\big(\Borel\cap (\boraxi\!\times\!\boraxi )\big)_\sigma$. It remains to check that 
$$\Borel\cap (\boraxi\!\times\!\boraxi )\! =\! (\Borel\cap\boraxi )\!\times\! (\Borel\cap\boraxi ).$$ 
The second set is clearly a subset of the first one. So assume that 
$R\! =\! A\!\times\! B\!\in\!\Borel\cap (\boraxi\!\times\!\boraxi )$. We may assume that $R$ is not empty. Then the projections $A$, $B$ are $\Ana$ since $R\!\in\!\Borel$. Lemma 2.2 gives 
$A',B'\!\in\!\Borel\cap\boraxi$ with $A\!\times\! B\!\subseteq\! A'\!\times\! B'\!\subseteq\! R\! =\! A\!\times\! B$. \hfill{$\square$}\bigskip

 Recall that if $A$ is a relation on $X$ and $D\!\subseteq\! X$, then $D$ is $A$-$discrete$ if 
$A\cap D^2\! =\!\emptyset$.\bigskip

\noindent\bf Proof of Theorem 2.1.\rm\ We apply Theorem 2.3 to $Y\! :=\! X$, $A_0\! :=\!\Delta (X)$ and 
$A_1\! :=\! A$. As 
$$(X,A)\preceq_{\borxi}\big(\omega ,\neg\Delta (\omega )\big)\mbox{,}$$ 
$\Delta (X)$ is separable from $A$ by a $(\boraxi\!\times\!\boraxi )_\sigma$ set. Theorem 2.3 gives 
$C_n,D_n\!\in\!\Borel\cap\boraxi$ such that $S\! :=\!\bigcup_{n\in\omega}~C_n\!\times\! D_n\!\in\!\Borel$ separates $\Delta (X)$ from $A$. As the set of codes for $\Borel\cap\boraxi$ subsets of $X$ is $\Ca$ (see Proposition 1.4 and Theorem 1.A in [Lo1]), the $\Borel$-selection theorem and the separation theorem imply that we may assume that the sequences $(C_n)$ and $(D_n)$ are $\Borel$. Note that 
$(C_n\cap D_n)$ is a $\Borel$ covering of $X$ into $A$-discrete $\Borel\cap\boraxi$ sets. As $X$ is $0$-dimensional we can reduce this covering into a $\Borel$ covering of $X$ consisting of $\Borel\cap\boraxi$ sets, which are in fact $\borxi$. This gives the desired partition.\hfill{$\square$}\bigskip

\noindent\bf Notation.\rm\ Following [Lo1], we define the following topologies on a $0$-dimensional recursively in $\alpha$ presented Polish space $X$, for any $\alpha\!\in\!\omega^\omega$. Let 
$T_1 (\alpha )$ be the topology of $X$, and, for $2\!\leq\!\xi\! <\!\omega_1$, $T_\xi (\alpha )$ be the topology generated by the $\Ana (\alpha )\cap\bormlxi$ subsets of $X$. The next proposition gives a reformulation of the inequality 
$(X,A)\preceq_{\Borel (\alpha )\cap\borxi}\big(\omega ,\neg\Delta (\omega )\big)$ of the Effective 
conjecture 1.

\begin{prop} Let $1\!\leq\!\xi\! <\!\omega_1^{\mbox{CK}}$, $X$ be a $0$-dimensional recursively presented Polish space, and $A$ be a $\Ana$ relation on $X$. Then 
$(X,A)\preceq_{\Borel\cap\borxi}\big(\omega ,\neg\Delta (\omega )\big)$ is equivalent to 
$\Delta (X)\cap\overline{A}^{T_\xi\times T_\xi}\! =\!\emptyset$.\end{prop}

\noindent\bf Proof.\rm\ Assume first that 
$(X,A)\preceq_{\Borel\cap\borxi}\big(\omega ,\neg\Delta (\omega )\big)$. Then there is a partition $(B_n)$ of $X$ into $A$-discrete $\Borel\cap\borxi$ sets. In particular, Theorem 1.A in [Lo1] implies that  $B_n$ is a countable union of $\Borel\cap\bormlxi$ sets if $\xi\!\geq\! 2$. In particular, $B_n$ is $T_\xi$-open and 
$\Delta (X)$ is disjoint from $\overline{A}^{T_\xi\times T_\xi}$ (even if $\xi\! =\! 1$).\bigskip

 Conversely, assume that $\Delta (X)\cap\overline{A}^{T_\xi\times T_\xi}\! =\!\emptyset$. Then each element $x$ of $X$ is contained in a $A$-discrete $\Ana\cap\bormlxi$ set (basic clopen set if 
$\xi\! =\! 1$). Lemma 2.2 implies that each element $x$ of $X$ is in fact contained in a $A$-discrete $\Borel\cap\bormlxi$ set if $\xi\!\geq\! 2$. It remains to apply Proposition 1.4 in [Lo1] and the $\Borel$-selection theorem to get the desired partition.\hfill{$\square$}\bigskip

 One can also hope for an effective strengthening of Conjecture 2 generalizing Theorem 1.5:\bigskip

\noindent\bf Effective conjecture 2\it\ Let $1\!\leq\!\xi\! <\!\omega_1$. Then there are\smallskip

\noindent - $0$-dimensional Polish spaces $\mathbb{X}^0_\xi ,\mathbb{X}^1_\xi$,\smallskip

\noindent - disjoint analytic subsets $\mathbb{A}^0_{\xi},\mathbb{A}^1_{\xi}$ of the space 
$\mathbb{X}^0_\xi\!\times\!\mathbb{X}^1_\xi$, not separable by a $(\boraxi\!\times\boraxi )_\sigma$ set,\smallskip

\noindent such that for any $\alpha\!\in\!\omega^\omega$ such that $1\!\leq\!\xi\! <\!\omega_1^\alpha$, for any recursively in $\alpha$ presented Polish spaces $X,Y$, and for any pair 
$A_0,A_1$ of disjoint $\Ana (\alpha )$ subsets of $X\!\times\! Y$, the following are equivalent:\smallskip

\noindent (a) The set $A_0$ cannot be separated from $A_1$ by a 
$(\boraxi\!\times\!\boraxi )_\sigma$ set.\smallskip

\noindent (b) The set $A_0$ cannot be separated from $A_1$ by a 
$\Borel (\alpha )\cap (\boraxi\!\times\!\boraxi )_\sigma$ set.\smallskip

\noindent (c) The set $A_0$ cannot be separated from $A_1$ by a 
$\boraone\big( T_\xi (\alpha )\!\times\! T_\xi (\alpha )\big)$ set.\smallskip

\noindent (d) $A_0\cap\overline{A_1}^{T_\xi (\alpha )\times T_\xi (\alpha )}\!\not=\!\emptyset$.\smallskip

\noindent (e) $(\mathbb{X}^0_\xi ,\mathbb{X}^1_\xi ,\mathbb{A}^0_\xi ,\mathbb{A}^1_\xi )\leq 
(X,Y,A_0,A_1)$.\rm

\vfill\eject

 In fact, the statements (a)-(d) are indeed equivalent:

\begin{thm} Let $1\!\leq\!\xi\! <\!\omega_1^{\mbox{CK}}$, $X, Y$ be recursively presented Polish spaces, and $A_0,A_1$ be disjoint $\Ana$ subsets of $X\!\times\! Y$. The following are equivalent:\smallskip

\noindent (a) The set $A_0$ cannot be separated from $A_1$ by a 
$(\boraxi\!\times\!\boraxi )_\sigma$ set.\smallskip

\noindent (b) The set $A_0$ cannot be separated from $A_1$ by a 
$\Borel\cap (\boraxi\!\times\!\boraxi )_\sigma$ set.\smallskip

\noindent (c) The set $A_0$ cannot be separated from $A_1$ by a 
$\boraone (T_\xi\!\times\! T_\xi )$ set.\smallskip

\noindent (d) $A_0\cap\overline{A_1}^{T_\xi\times T_\xi}\!\not=\!\emptyset$.\end{thm}

\noindent\bf Proof.\rm\ Theorem 2.3 implies that (a) is indeed equivalent to (b). It also implies, using the proof of Proposition 2.4, that (c) implies (a), and the converse is clear. It is also clear that (c) and (d) are equivalent.\hfill{$\square$}\bigskip

 A consequence of this is that Conjecture 2 and the Effective conjecture 2 are equivalent.

\section{$\!\!\!\!\!\!$ The case $\xi\! =\! 1$}

\bf (A) Continuous colorings\rm\bigskip

 As in [L3], we can separate Conjecture 1 in two parts. We introduce the following notion, that will help us to characterize the relations $\mathbb{A}$ for which there is a continuous homomorphism from $\mathbb{A}$ into any relation without countable continuous coloring:

\begin{defi} Let $\xi$ be a countable ordinal, ${\bf\Pi}^0_0\! :=\!\borone$, and $\mathbb{X}$ be a $0$-dimensional Polish  space. A family $\cal F$ of subsets of $\mathbb{X}$ is $\xi$-$disjoint$ if the elements of $\cal F$ are $\bormxi$ and pairwise disjoint.\end{defi}

 The first part ensures the existence of complicated examples.
 
\begin{lem} (a) Assume that $(C^\varepsilon_i)_{(\varepsilon ,i)\in 2\times\omega}$ is a $0$-disjoint family of subsets of the space $\mathbb{X}$ such that 
$\mathbb{X}\!\setminus\! (\bigcup_{(\varepsilon ,i)\in 2\times\omega}~C^\varepsilon_i)\!\not=\!\emptyset$  and no clopen set meeting 
$\mathbb{X}\!\setminus\! (\bigcup_{(\varepsilon ,i)\in 2\times\omega}~C^\varepsilon_i)$ is 
$(\bigcup_{i\in\omega}~C^0_i\!\times\! C^1_i)$-discrete. Then 
$(\mathbb{X},\bigcup_{i\in\omega}~C^0_i\!\times\! C^1_i)\not\preceq_{\borone}
\big(\omega ,\neg\Delta (\omega )\big)$.\smallskip

\noindent (b) There is a $0$-disjoint family $(C^\varepsilon_i)_{(\varepsilon ,i)\in 2\times\omega}$ of subsets of $2^\omega$ satisfying the assumption (and thus the conclusion) of (a).\end{lem}

\noindent\bf Proof.\rm\ (a) We argue by contradiction, which gives $f\! :\!\mathbb{X}\!\rightarrow\!\omega$ continuous such that $f(x)\!\not=\! f(y)$ if $(x,y)\!\in\!\bigcup_{i\in\omega}~C^0_i\!\times\! C^1_i$. We set $D_k\! :=\! f^{-1}(\{ k\})$, so that $(D_k)_{k\in\omega}$ is a partition of $\mathbb{X}$ into clopen sets discrete for $\bigcup_{i\in\omega}~C^0_i\!\times\! C^1_i$. Choose 
$z\!\in\!\mathbb{X}\!\setminus\! (\bigcup_{(\varepsilon ,i)\in 2\times\omega}~C^\varepsilon_i)$, and $k$ with $z\!\in\! D_k$. This gives $(x,y)\!\in\! (\bigcup_{i\in\omega}~C^0_i\!\times\! C^1_i)\cap D_k^2$, which is absurd.\bigskip

\noindent (b) We set $C^\varepsilon_i\! :=\! N_{0^{2i+\varepsilon}1}$, so that 
$\bigcup_{i\in\omega}~C^0_i\!\times\! C^1_i\! =\!\{ (0^{2i}1\alpha ,0^{2i+1}1\beta )\mid 
i\!\in\!\omega\mbox{ and }\alpha ,\beta\!\in\! 2^\omega\}$. Note that 
$\{ 0^\infty\}\! =\!\mathbb{X}\!\setminus\! (\bigcup_{(\varepsilon ,i)\in 2\times\omega}~C^\varepsilon_i)$. If 
$C$ is a clopen neighborhood of $0^\infty$, then $N_{0^i}\!\subseteq\! C$ if $i$ is big enough. This gives an integer $i$ with $(0^{2i}1^\infty ,0^{2i+1}1^\infty )\!\in\! (\bigcup_{i\in\omega}~C^0_i\!\times\! C^1_i)\cap C^2$.\hfill{$\square$}\bigskip
 
 The second part ensures the existence of the continuous homomorphism.
 
\begin{lem} Let $\mathbb{X}$ be a $0$-dimensional Polish  space, 
$(C^\varepsilon_i)_{(\varepsilon ,i)\in 2\times\omega}$ be a $0$-disjoint family of subsets of 
$\mathbb{X}$, $X$ be a $0$-dimensional Polish space, and $A$ be a relation on $X$. Then one of the following holds:\smallskip  

\noindent (a) $(X,A)\preceq_{\borone}\big(\omega ,\neg\Delta (\omega )\big)$,\smallskip  

\noindent (b) $(\mathbb{X},\bigcup_{i\in\omega}~C^0_i\!\times\! C^1_i)\preceq_{\boraone}(X,A)$.
\end{lem}

\noindent\bf Proof.\rm\ Assume that (a) does not hold. Let us fix a compatible complete metric on $X$. In the sequel, the diameter will refer to this metric (this will also be the case in all the proofs where diameters are involved to come). We enumerate a basis $\big( N(p,X)\big)_{p\in\omega}$ for the topology of $X$ made of clopen sets.\bigskip

\noindent $\bullet$ We build\medskip
 
\noindent - an increasing sequence of integers $(p_i)_{i\in\omega}$,\smallskip
 
\noindent - a sequence $(x_k)_{k\in\omega}$ of points of $X$.\bigskip
 
 We want these objects to satisfy the following conditions:
$$\begin{array}{ll}
& (1)~(x_{2i},x_{2i+1})\!\in\! A\cap N(p_i,X)^2\cr
& (2)~N(p_{i+1},X)\!\subseteq\! N(p_i,X)\cr
& (3)~\mbox{diam}\big( N(p_i,X)\big)\!\leq\! 2^{-i}\cr
& (4)~\mbox{There is no covering of }N(p_i,X)\mbox{ consisting of }A\mbox{-discrete clopen subsets of }X
\end{array}$$ 
$\bullet$ Assume that this is done. Then we can define a point $x$ of $X$ by 
$\{ x\}\! =\!\bigcap_{i\in\omega}~N(p_i,X)$. Note that $(x_k)_{k\in\omega}$ tends to $x$. We define 
$f\! :\!\mathbb{X}\!\rightarrow\! X$ by $f(z)\! :=\! x$ if 
$z\!\notin\!\bigcup_{(\varepsilon ,i)\in 2\times\omega}~C^\varepsilon_i$, $f(z)\! :=\! x_{2i+\varepsilon}$ if 
$z\!\in\! C^\varepsilon_i$. Note that $f$ is continuous. Moreover, 
$\big( f(y),f(z)\big)\! =\! (x_{2i},x_{2i+1})\!\in\! A$ if $(y,z)\!\in\! C^0_i\!\times\! C^1_i$, so that (b) holds.\bigskip

\noindent $\bullet$ Let us prove that the construction is possible. We set $N(p_{-1},X)\! :=\! X$. Assume that $(p_i)_{i<l}$ and $(x_{2i},x_{2i+1})_{i<l}$ satisfying (1)-(4) have been constructed, which is the case for $l\! =\! 0$. We choose a covering of $N(p_{l-1},X)$ with basic clopen sets of diameter at most $2^{-l}$, contained in $N(p_{l-1},X)$. Then one of these basic sets, say $N(p_l,X)$, satisfies (4). It remains to choose $(x_{2l},x_{2l+1})$ in the set $A\cap N(p_l,X)^2$.\hfill{$\square$}\bigskip

 We set $\mathbb{X}_1\! :=\! 2^\omega$ and $\mathbb{A}_1\! :=\!\{ (0^{2i}1\alpha ,0^{2i+1}1\beta )\mid 
i\!\in\!\omega\mbox{ and }\alpha ,\beta\!\in\! 2^\omega\}\! =\!\bigcup_{i\in\omega}~C^0_i\!\times\! C^1_i$, so that $\mathbb{A}_1$ is a $\boraone$ relation on $\mathbb{X}_1$.

\begin{cor} Let $X$ be a $0$-dimensional Polish space, and $A$ be a relation on $X$. Then exactly one of the following holds:\smallskip  

\noindent (a) $(X,A)\preceq_{\borone}\big(\omega ,\neg\Delta (\omega )\big)$,\smallskip  

\noindent (b) $(\mathbb{X}_1,\mathbb{A}_1)\preceq_{\boraone}(X,A)$.\smallskip

 Moreover, there are a non $0$-dimensional Polish space $X$, and a closed relation $A$ on $X$, for which neither (a), nor (b) holds (with this couple $(\mathbb{X}_1,\mathbb{A}_1)$ or any other). There are also a $0$-dimensional Polish space $X$, and a relation $A$ on $X$ (a difference of two closed sets), for which it is not possible to have $f$ one-to-one in (b) (with this couple $(\mathbb{X}_1,\mathbb{A}_1)$ or any other).\end{cor}
 
\noindent\bf Proof.\rm\ Note first that (a) and (b) cannot hold simultaneously, by Lemma 3.2. Lemma 3.3 implies that (a) or (b) holds.\bigskip

\noindent $\bullet$ Consider now $X\! :=\!\mathbb{R}$ and $A\! :=\!\{ (0,1)\}$. Then (a) does not hold since 
$\mathbb{R}$ is connected. If (b) holds, then we must have $f(0^{2i}1\alpha )\! =\! 0$ and 
$f(0^{2i+1}1\beta )\! =\! 1$. By continuity of $f$, we get $f(0^\infty )\! =\! 0\! =\! 1$.\bigskip

 This would be the same with any $(\mathbb{X}_1,\mathbb{A}_1)$. Indeed, as 
$(\mathbb{X}_1,\mathbb{A}_1)\not\preceq_{\borone}\big(\omega ,\neg\Delta (\omega )\big)$, we have 
$\overline{\Pi_0[\mathbb{A}_1]}\cap\overline{\Pi_1[\mathbb{A}_1]}\!\not=\!\emptyset$, since otherwise there would be a clopen subset $C$ of the $0$-dimensional space $\mathbb{X}_1$ separating 
$\overline{\Pi_0[\mathbb{A}_1]}$ from $\overline{\Pi_1[\mathbb{A}_1]}$, and we would have 
$\Delta (\mathbb{X}_1)\!\subseteq\! C^2\cup (\neg C)^2\!\subseteq\!\neg\mathbb{A}_1$. So we can choose ${x\!\in\!\overline{\Pi_0[\mathbb{A}_1]}\cap\overline{\Pi_1[\mathbb{A}_1]}}$, 
$x_{2i}\!\in\!\Pi_0[\mathbb{A}_1]$ such that $(x_{2i})$ tends to $x$, 
$y_{2i+1}\!\in\!\Pi_1[\mathbb{A}_1]$ such that $(y_{2i+1})$ tends to $x$, 
$y_{2i}$ with $(x_{2i},y_{2i})\!\in\!\mathbb{A}_1$, and 
$x_{2i+1}$ with $(x_{2i+1},y_{2i+1})\!\in\!\mathbb{A}_1$. Then $f(x_{2i})\! =\! 0$, $f(y_{2i+1})\! =\! 1$ and we conclude as before.\bigskip

\noindent $\bullet$ Consider $X\! :=\! 2^\omega$ and 
$A\! :=\!\{ 0^\infty\}\!\times\! (2^\omega\!\setminus\!\{ 0^\infty\} )$. Then (a) does not hold since if a clopen subset $C$ of $2^\omega$ contains $0^\infty$, then it contains also some $\alpha\!\not=\! 0^\infty$, so that 
$(0^\infty ,\alpha )\!\in\! A\cap C^2$. If (b) holds, then $f(0^{2i}1\alpha )\! =\! 0^\infty$ for each integer $i$ and $f$ is not one-to-one.\bigskip

 This argument works as soon as $\Pi_0[\mathbb{A}_1]$ has at least two elements. If we argue in the other factor, then we see that an example $(\mathbb{X}_1,\mathbb{A}_1)$ with injectivity must satisfy that $\mathbb{A}_1$ is a singleton $\{ (\alpha ,\beta )\}$. As 
$(\mathbb{X}_1,\mathbb{A}_1)\preceq_{\boraone}(2^\omega ,\mathbb{G}_0)$, $\alpha\!\not=\!\beta$. So take a clopen subset $C$ of $\mathbb{X}_1$ containing $\alpha$ but not $\beta$. Then 
$\Delta (\mathbb{X}_1)\!\subseteq\! C^2\cup (\neg C)^2\!\subseteq\!\neg\mathbb{A}_1$.\hfill{$\square$}\bigskip

 The notion of a $0$-disjoint family is essential in the following sense:
  
\begin{prop} Let $\mathbb{X}$ be a $0$-dimensional Polish space, and $\mathbb{A}$ be a relation on 
$\mathbb{X}$. The following are equivalent:\smallskip

\noindent (a) For any $0$-dimensional Polish space $X$, and any relation $A$ on $X$,
$$(X,A)\not\preceq_{\borone}\big(\omega ,\neg\Delta (\omega )\big)\Rightarrow 
(\mathbb{X},\mathbb{A})\preceq_{\boraone}(X,A).$$
(b) There is a $0$-disjoint family $(C^\varepsilon_i)_{(\varepsilon ,i)\in 2\times\omega}$ of subsets of 
$\mathbb{X}$ such that $\mathbb{A}\!\subseteq\!\bigcup_{i\in\omega}~C^0_i\!\times\! C^1_i$.
\end{prop}

\noindent\bf Proof.\rm\ (a) $\Rightarrow$ (b) We set $X\! :=\!\mathbb{X}_1$ and $A\! :=\!\mathbb{A}_1$. By Lemma 3.2, we get $f\! :\!\mathbb{X}\!\rightarrow\! 2^\omega$ such that 
$\mathbb{A}\!\subseteq\! (f\!\times\!f)^{-1}(\mathbb{A}_1)$. We set 
$C^\varepsilon_i\! :=\! f^{-1}(N_{0^{2i+\varepsilon}1})$.\bigskip

\noindent (b) $\Rightarrow$ (a) By Lemma 3.3 we get 
$(\mathbb{X},\bigcup_{i\in\omega}~C^0_i\!\times\! C^1_i)\preceq_{\boraone}(X,A)$, so that 
$(\mathbb{X},\mathbb{A})\preceq_{\boraone}(X,A)$.\hfill{$\square$}\bigskip

\noindent\bf (B) Countable unions of open rectangles (i.e., open sets)\rm\bigskip

 The content here is completely trivial. It is just the fact that a subset of a metric space is not open exactly when it contains a point that we can approximate by a countable sequence contained in its complement. We give some statements since the situation will be more complicated in the case $\xi\! =\! 2$. As in (A) we can characterize the tuples $(\mathbb{X}^0,\mathbb{X}^1,\mathbb{A}^0,\mathbb{A}^1)$ $\leq$-below any tuple $(X,Y,A_0,A_1)$ with $A_0$ not separable from $A_1$ by a $(\boraone\!\times\!\boraone )_\sigma$ set.
  
\begin{lem} (a) Assume that $(C^\varepsilon_i)_{i\in\omega}$ is a $0$-disjoint family of subsets of the space $\mathbb{X}^\varepsilon$ such that 
$\big(\mathbb{X}^0\!\setminus\! (\bigcup_{i\in\omega}~C^0_i)\big)\!\times\!
\big(\mathbb{X}^1\!\setminus\! (\bigcup_{i\in\omega}~C^1_i)\big)\!\not=\!\emptyset$ and no open set meeting $\big(\mathbb{X}^0\!\setminus\! (\bigcup_{i\in\omega}~C^0_i)\big)\!\times\!
\big(\mathbb{X}^1\!\setminus\! (\bigcup_{i\in\omega}~C^1_i)\big)$ is disjoint from 
$\bigcup_{i\in\omega}~C^0_i\!\times\! C^1_i$. Then 
$\big(\mathbb{X}^0\!\setminus\! (\bigcup_{i\in\omega}~C^0_i)\big)\!\times\!
\big(\mathbb{X}^1\!\setminus\! (\bigcup_{i\in\omega}~C^1_i)\big)$ is not separable from 
$\bigcup_{i\in\omega}~C^0_i\!\times\! C^1_i$ by a $(\boraone\!\times\!\boraone )_\sigma$ set.\smallskip

\noindent (b) There are $0$-disjoint families of subsets of $2^\omega$ satisfying the assumption (and thus the conclusion) of (a).\end{lem}
 
\noindent\bf Proof.\rm\ (a) is obvious.\bigskip

\noindent (b) We set $C^\varepsilon_i\! :=\! N_{0^i1}$, so that 
$\bigcup_{i\in\omega}~C^0_i\!\times\! C^1_i\! =\!\{ (0^i1\alpha ,0^i1\beta )\mid 
i\!\in\!\omega\mbox{ and }\alpha ,\beta\!\in\! 2^\omega\}$. Note that 
$\{ 0^\infty\}\! =\!\mathbb{X}^\varepsilon\!\setminus\! (\bigcup_{i\in\omega}~C^\varepsilon_i)$. If 
$O$ is an open neighborhood of $(0^\infty ,0^\infty )$, then $N^2_{0^i}\!\subseteq\! O$ if $i$ is big enough. This gives an integer $i$ with 
$(0^i1^\infty ,0^i1^\infty )\!\in\! (\bigcup_{i\in\omega}~C^0_i\!\times\! C^1_i)\cap O$.\hfill{$\square$}
 
\begin{lem} Let $\mathbb{X}^0,\mathbb{X}^1$ be $0$-dimensional Polish  spaces, 
$(C^\varepsilon_i)_{i\in\omega}$ be a $0$-disjoint family of subsets of $\mathbb{X}^\varepsilon$, 
$X,Y$ be Polish spaces, and $A_0,A_1$ be disjoint subsets of $X\!\times\! Y$. Then one of the following holds:\smallskip  

\noindent (a) $A_0$ is separable from $A_1$ by a $(\boraone\!\times\!\boraone )_\sigma$ set,\smallskip  

\noindent (b) $\big(\mathbb{X}^0,\mathbb{X}^1,
\big(\mathbb{X}^0\!\setminus\! (\bigcup_{i\in\omega}~C^0_i)\big)\!\times\!
\big(\mathbb{X}^1\!\setminus\! (\bigcup_{i\in\omega}~C^1_i)\big),
\bigcup_{i\in\omega}~C^0_i\!\times\! C^1_i\big)\leq (X,Y,A_0,A_1)$.
\end{lem}

\noindent\bf Proof.\rm\ Assume that (a) does not hold. Pick $(x,y)\!\in\! A_0\cap\overline{A_1}$, and 
$(x_i,y_i)$ in $A_1$ tending to $(x,y)$. We define $f\! :\!\mathbb{X}^0\!\rightarrow\! X$ by $f(z)\! :=\! x$ if 
$z\!\notin\!\bigcup_{i\in\omega}~C^0_i$, $x_i$ if $z\!\in\! C^0_i$. Note that $f$ is continuous. Similarly, we define $g\! :\!\mathbb{X}^1\!\rightarrow\! Y$, so that (b) holds.\hfill{$\square$}\bigskip

 We define $\mathbb{X}^\varepsilon_1\! :=\! 2^\omega$, 
$\mathbb{A}^0_1\! :=\!\{ (0^\infty ,0^\infty )\}\! =\! 
\big(\mathbb{X}^0_1\!\setminus\! (\bigcup_{i\in\omega}~C^0_i)\big)\!\times\!
\big(\mathbb{X}^1_1\!\setminus\! (\bigcup_{i\in\omega}~C^1_i)\big)$ and also 
$\mathbb{A}^1_1\! :=\!\{ (0^i1\alpha ,0^i1\beta )\mid i\!\in\!\omega\mbox{ and }
\alpha ,\beta\!\in\! 2^\omega\}\! =\!\bigcup_{i\in\omega}~C^0_i\!\times\! C^1_i$. As in (A) we get the two following consequences:

\begin{cor} Let $X,Y$ be Polish spaces, and $A_0,A_1$ be disjoint subsets of $X\!\times\! Y$. Then exactly one of the following holds:\smallskip  

\noindent (a) $A_0$ is separable from $A_1$ by a $(\boraone\!\times\!\boraone )_\sigma$ set,\smallskip  

\noindent (b) $(\mathbb{X}^0_1,\mathbb{X}^1_1,\mathbb{A}^0_1,\mathbb{A}^1_1)\leq (X,Y,A_0,A_1)$.
\end{cor}
  
\begin{prop} Let $\mathbb{X}^0,\mathbb{X}^1$ be $0$-dimensional Polish spaces, and 
$\mathbb{A}^0,\mathbb{A}^1\!\subseteq\!\mathbb{X}^0\!\times\!\mathbb{X}^1$ be disjoint. The following are equivalent:\smallskip

\noindent (a) For any Polish spaces $X,Y$, and any $A_0,A_1\!\subseteq\! X\!\times\! Y$ disjoint,
$$A_0\mbox{ is not separable from }A_1\mbox{ by a }(\boraone\!\times\!\boraone )_\sigma\mbox{ set }\Rightarrow 
(\mathbb{X}^0,\mathbb{X}^1,\mathbb{A}^0,\mathbb{A}^1)\leq (X,Y,A_0,A_1).$$
(b) There is a $0$-disjoint family $(C^\varepsilon_i)_{i\in\omega}$ of subsets of 
$\mathbb{X}^\varepsilon$ such that the inclusions 
$\mathbb{A}^1\!\subseteq\!\bigcup_{i\in\omega}~C^0_i\!\times\! C^1_i$ and $\mathbb{A}^0\!\subseteq\!
\big(\mathbb{X}^0\!\setminus\! (\bigcup_{i\in\omega}~C^0_i)\big)\!\times\!
\big(\mathbb{X}^1\!\setminus\! (\bigcup_{i\in\omega}~C^1_i)\big)$ hold.
\end{prop}

\section{$\!\!\!\!\!\!$ The case $\xi\! =\! 2$}

\bf (A) Baire class one colorings\rm

\begin{lem} (a) Assume that $(C^\varepsilon_i)_{(\varepsilon ,i)\in 2\times\omega}$ is a $1$-disjoint family of subsets of $\mathbb{X}$ such that no non-empty clopen subset of $\mathbb{X}$ is 
$(\bigcup_{i\in\omega}~C^0_i\!\times\! C^1_i)$-discrete. Then 
$(\mathbb{X},\bigcup_{i\in\omega}~C^0_i\!\times\! C^1_i)\preceq_{\bortwo}
\big(\omega ,\neg\Delta (\omega )\big)$.\smallskip

\noindent (b) There is a $1$-disjoint family $(C^\varepsilon_i)_{(\varepsilon ,i)\in 2\times\omega}$ of subsets of $\omega^\omega$ satisfying the assumption (and thus the conclusion) of (a).\end{lem}

\noindent\bf Proof.\rm\ (a) We argue by contradiction, which gives a $\bortwo$-measurable map  
$f\! :\!\mathbb{X}\!\rightarrow\!\omega$ with $f(x)\!\not=\! f(y)$ if 
$(x,y)\!\in\!\bigcup_{i\in\omega}~C^0_i\!\times\! C^1_i$. We set $D_k\! :=\! f^{-1}(\{ k\})$, so that 
$(D_k)_{k\in\omega}$ is a partition of $\mathbb{X}$ into $\bortwo$ sets discrete for 
$\bigcup_{i\in\omega}~C^0_i\!\times\! C^1_i$. By Baire's theorem, there are an integer $k$ and a nonempty clopen subset $C$ of $\mathbb{X}$ such that $D_k$ contains $C$. This gives 
$(x,y)\!\in\! (\bigcup_{i\in\omega}~C^0_i\!\times\! C^1_i)\cap C^2\!\subseteq\! 
(\bigcup_{i\in\omega}~C^0_i\!\times\! C^1_i)\cap D_k^2$, which is absurd.\bigskip

\noindent (b) Let $b\! :\!\omega\!\rightarrow\!\omega^{<\omega}$ be a bijection. We set 
$C^\varepsilon_i\! :=\!\big\{ b(i)\big( 2\vert b(i)\vert\! +\!\varepsilon\big)^\infty\big\}$, so that 
$$\bigcup_{i\in\omega}~C^0_i\!\times\! C^1_i\! =\!
\big\{\big( u(2\vert u\vert )^\infty ,u(2\vert u\vert\! +\! 1)^\infty\big)\mid u\!\in\!\omega^{<\omega}\big\} .$$ 
If $\emptyset\!\not=\! C\!\in\!\borone (\omega^\omega )$, then $C$ contains some basic clopen set $N_u$, and $\big( u(2\vert u\vert )^\infty ,u(2\vert u\vert\! +\! 1)^\infty\big)$ is in 
$(\bigcup_{i\in\omega}~C^0_i\!\times\! C^1_i)\cap C^2$.\hfill{$\square$}\bigskip

\noindent\bf Remark.\rm\ There are a $1$-disjoint family 
$(C^\varepsilon_i)_{(\varepsilon ,i)\in 2\times\omega}$ of subsets of $\omega^\omega$, and a relation $A$ on $\omega^\omega$ such that 
$(\omega^\omega ,A)\not\preceq_{\bortwo}\big(\omega ,\neg\Delta (\omega )\big)$ and 
$(\omega^\omega ,\bigcup_{i\in\omega}~C^0_i\!\times\! C^1_i)\not\preceq_{\boraone}
(\omega^\omega ,A)$, so that Lemma 3.3 cannot be extended to $\boratwo$-measurable countable colorings.\bigskip

 Indeed, we set $C^\varepsilon_i\! :=\!\{ u(2i\! +\!\varepsilon )^\infty\mid u\!\in\!\omega^i\}$ and 
$A\! :=\!\big\{\big( u(2\vert u\vert )^\infty ,u(2\vert u\vert\! +\! 1)^\infty\big)\mid u\!\in\!\omega^{<\omega}\big\}$. Then $(C^\varepsilon_i)_{(\varepsilon ,i)\in 2\times\omega}$ is clearly a $1$-disjoint family. Lemma 4.1 gives the first assertion. For the second assertion, assume, towards a contradiction, that 
$f\! :\!\omega^\omega\!\rightarrow\!\omega^\omega$ is continuous and satisfies the inclusion 
$\bigcup_{i\in\omega}~C^0_i\!\times\! C^1_i\!\subseteq\! (f\!\times\! f)^{-1}(A)$. If $i\!\in\!\omega$, then there is $u_i\!\in\!\omega^{<\omega}$ with 
$$f[C^0_i]\!\times\! f[C^1_i]\!\subseteq\!
\big\{\big( u_i(2\vert u_i\vert )^\infty ,u_i(2\vert u_i\vert\! +\! 1)^\infty\big)\big\}.$$
In particular, for any $\alpha ,\beta\!\in\!\omega^\omega$ we get 
$$\begin{array}{ll}
\big( f(\alpha ),f(\beta )\big)\!\!\!\!
& \! =\!\mbox{lim}_{i\rightarrow\infty}~
\Big( f\big( (\alpha\vert i)(2i)^\infty\big),f\big( (\beta\vert i)(2i\! +\! 1)^\infty\big)\Big)\cr\cr
& \! =\!\mbox{lim}_{i\rightarrow\infty}~
\big( u_i(2\vert u_i\vert )^\infty ,u_i(2\vert u_i\vert\! +\! 1)^\infty\big).
\end{array}$$ 
But this implies that $f$ is constant, which is absurd. To fix  this, we refine the notion of a $\xi$-disjoint family.

\begin{defi} Let $1\!\leq\!\xi\! <\!\omega_1$. A $\xi$-disjoint family 
$(C^\varepsilon_i)_{(\varepsilon ,i)\in 2\times\omega}$ of subsets of a $0$-dimensional Polish space 
$\mathbb{X}$ is said to be $comparing$ if for each integer $q$ there is a partition 
$(O_q^p)_{p\in\omega}$ of $\mathbb{X}$ into $\borxi$ sets such that, for each $i\!\in\!\omega$,\smallskip

\noindent (a) if $q\! <\! i$, then there is $p_q^i\!\in\!\omega$ such that 
$C^0_i\cup C^1_i\!\subseteq\! O_q^{p_q^i}$,\smallskip

\noindent (b) if $q\!\geq\! i$ and $\varepsilon\!\in\! 2$, then 
$C^\varepsilon_i\!\subseteq\! O_q^{2i+\varepsilon}$.\end{defi}

\begin{lem} There is a comparing $1$-disjoint family 
$(C^\varepsilon_i)_{(\varepsilon ,i)\in 2\times\omega}$ of subsets of $\omega^\omega$ satisfying the assumption (and thus the conclusion) of Lemma 4.1.(a).\end{lem}

\noindent\bf Proof.\rm\ Let $b\! :\!\omega\!\rightarrow\!\omega^{<\omega}$ be a bijection satisfying  
$b^{-1}(s)\!\leq\!b^{-1}(t)$ if $s\!\subseteq\! t$. It can be built as follows. Let $(p_q)_{q\in\omega}$ be the sequence of prime numbers, and $I\! :\!\omega^{<\omega}\!\rightarrow\!\omega$ defined by 
$I(s)\! :=\! p_0^{s(0)+1}...p_{\vert s\vert -1}^{s(\vert s\vert -1)+1}$ if $s\!\not=\!\emptyset$, and 
$I(\emptyset )\! :=\! 1$. Note that $I$ is one-to-one, so that there is an increasing bijection 
$\varphi\! :\! I[\omega^{<\omega}]\!\rightarrow\!\omega$. We set 
$b\! :=\! (\varphi\circ I)^{-1}\! :\!\omega\!\rightarrow\!\omega^{<\omega}$. We define 
$(C^\varepsilon_i)_{(\varepsilon ,i)\in 2\times\omega}$ as in the proof of Lemma 4.1.(b), so that 
$(C^\varepsilon_i)_{(\varepsilon ,i)\in 2\times\omega}$ is a $1$-disjoint family. It remains to see that 
$(C^\varepsilon_i)_{(\varepsilon ,i)\in 2\times\omega}$ is comparing. We set
$$O^p_q\! :=\!\left\{\!\!\!\!\!\!\!
\begin{array}{ll}
& N_{b(i)(2\vert b(i)\vert+\varepsilon )^{\mbox{max}_{l\leq q}~(\vert b(l)\vert +1)-\vert b(i)\vert}}
\mbox{ if }p\! =\! 2i\! +\!\varepsilon\!\leq\! 2q\! +\! 1\mbox{,}\cr\cr
& \omega^\omega\!\setminus\! (\bigcup_{p'\leq 2q+1}~O^{p'}_q)\mbox{ if }p\! =\! 2q\! +\! 2\mbox{,}\cr\cr
& \emptyset\mbox{ if }p\!\geq\! 2q\! +\! 3\mbox{,}
\end{array}
\right.$$
so that $(O_q^p)_{p\in\omega}$ is a partition of $\omega^\omega$ into $\borone$ sets. Note that (b) is fulfilled. If $q\! <\! i$, then there is at most one couple 
$(j,\varepsilon )\!\in\! (q\! +\! 1)\!\times\! 2$ such that 
$b(j)(2\vert b(j)\vert+\varepsilon )^{\mbox{max}_{l\leq q}~(\vert b(l)\vert +1)-\vert b(j)\vert}$ is compatible with $b(i)$. If it exists and if $\vert b(i)\vert\!\geq\!\mbox{max}_{l\leq q}~(\vert b(l)\vert +1)$, then 
$C^0_i\cup C^1_i\!\subseteq\! O_q^{2j+\varepsilon}$ and we set $p_q^i\! :=\! 2j\! +\!\varepsilon$. Otherwise, $C^0_i\cup C^1_i\!\subseteq\! O_q^{2q+2}$ and we set $p_q^i\! :=\! 2q\! +\! 2$.
\hfill{$\square$}\bigskip

 We have a stronger result than Conjecture 1, in the sense that we do not need any regularity assumption on $A$, neither that $X$ is $0$-dimensional.
 
\begin{lem} Let $\mathbb{X}$ be a $0$-dimensional Polish space, 
$(C^\varepsilon_i)_{(\varepsilon ,i)\in 2\times\omega}$ be a comparing $1$-disjoint family of subsets of 
$\mathbb{X}$, $X$ be a Polish space, and $A$ be a relation on $X$. Then one of the following holds:\smallskip  

\noindent (a) $(X,A)\preceq_{\bortwo}\big(\omega ,\neg\Delta (\omega )\big)$,\smallskip  

\noindent (b) $(\mathbb{X},\bigcup_{i\in\omega}~C^0_i\!\times\! C^1_i)\preceq_{\boraone}(X,A)$.
\end{lem}

\noindent\bf Proof.\rm\ If $A$ is not a digraph, then choose $x_0$ with $(x_0,x_0)\!\in\! A$, and put 
$f(x)\! :=\! x_0$. So we may assume that $A$ is a digraph. We set 
$$U\! :=\!\bigcup\Big\{ V\!\in\!\boraone (X)\mid\exists D\!\in\!\bormone (\omega\!\times\! X)~~
V\!\subseteq\!\bigcup_{p\in\omega}~D_p\mbox{ and }\forall p\!\in\!\omega ~~
A\cap D^2_p\! =\!\emptyset\Big\}.$$
\bf Case 1.\rm\ $U\! =\! X$.\bigskip

 There is a countable covering of $X$ into $A$-discrete $\boratwo$ sets. We just have to reduce them to get a partition showing that (a) holds.\bigskip

\noindent\bf Case 2.\rm\ $U\!\not=\! X$.\bigskip

 Then $Y\! :=\! X\!\setminus\! U$ is a nonempty closed subset of $X$.\bigskip

\noindent\bf Claim\it\ If $V\!\in\!\boraone (X)$ meets $Y$, then $V\cap Y$ is not $A$-discrete.\rm\bigskip
 
 We argue by contradiction. As $V\cap U$ can be covered with some $\bigcup_{p\in\omega}~D_p$'s, so is $V$. Thus $V\!\subseteq\! U$, so that $V\cap Y\!\subseteq\! U\!\setminus\! U\! =\!\emptyset$, which is the desired contradiction.\hfill{$\diamond$}

\vfill\eject

\noindent $\bullet$ We construct a family $(x_u)_{u\in\omega^{<\omega}}$ of points of $Y$, and a family 
$(X_u)_{u\in\omega^{<\omega}}$ of open subsets of $Y$. We want these objects to satisfy the following conditions:
$$\begin{array}{ll}
& (1)~x_u\!\in\! X_u\cr
& (2)~\overline{X_{up}}\!\subseteq\! X_u\cr
& (3)~\mbox{diam}(X_u)\!\leq\! 2^{-\vert u\vert}\cr
& (4)~(x_{u(2n)},x_{u(2n+1)})\!\in\! A\mbox{ if }u\!\in\!\omega^n\cr
& (5)~x_{u(2n+\varepsilon )}\! =\! x_u\mbox{ if }u\!\notin\!\omega^n\mbox{ and }\varepsilon\!\in\! 2
\end{array}$$ 
$\bullet$ Assume that this is done. We define $f\! :\!\mathbb{X}\!\rightarrow\! Y\!\subseteq\! X$ by 
$$\{ f(x)\}\! :=\!\bigcap_{q\in\omega}~\overline{X_{p_0...p_{q-1}}}\! =\!
\bigcap_{q\in\omega}~X_{p_0...p_{q-1}}\mbox{,}$$
where $p_i$ satisfies $x\!\in\! O^{p_i}_i$ witnessing comparability, so that $f$ is continuous. Note that 
$f(x)$ is the limit of $x_{p_0...p_{q-1}}$, and that 
$x_{u(2\vert u\vert +\varepsilon )}\! =\! x_{u(2\vert u\vert +\varepsilon )^2}\! =\! ...\! =\! 
x_{u(2\vert u\vert +\varepsilon )^{q+1}}$ 
for each $(u,\varepsilon )\!\in\!\omega^{<\omega}\!\times\! 2$. Thus 
$f(x)\! =\!\mbox{lim}_{q\rightarrow\infty}~x_{u(2\vert u\vert +\varepsilon )^{q+1}}\! =\! 
x_{u(2\vert u\vert +\varepsilon )}$ if $x\!\in\! C^\varepsilon_i$ and $u\! :=\! p_0^i...p_{i-1}^i$, and 
$$\big( f(x),f(y)\big)\! =\! (x_{u(2\vert u\vert )},x_{u(2\vert u\vert +1)})\!\in\! A$$ 
if $(x,y)\!\in\! C^0_i\!\times\! C^1_i$. So (b) holds.\bigskip

\noindent $\bullet$ Let us prove that the construction is possible. We choose $x_\emptyset\!\in\! Y$ and 
an open neighborhood $X_\emptyset$ of $x_\emptyset$ in $Y$, of diameter at most $1$. Assume that 
$(x_u)_{u\in\omega^{\leq l}}$ and $(X_u)_{u\in\omega^{\leq l}}$ satisfying (1)-(5) have been constructed, which is the case for $l\! =\! 0$.\bigskip

 An application of the claim gives $(x_{u(2l)},x_{u(2l+1)})\!\in\! A\cap X_u^2$ if $u\!\in\!\omega^l$. We satisfy (5), so that the definition of the $x_u$'s is complete. Note that  $x_u\!\in\! X_{u\vert l}$ if 
$u\!\in\!\omega^{l+1}$.\bigskip

 We choose an open neighborhood $X_u$ of $x_u$ in $Y$, of diameter at most $2^{-l-1}$, ensuring the inclusion $\overline{X_u}\!\subseteq\! X_{u\vert l}$.\hfill{$\square$}\bigskip
 
 We set $\mathbb{X}_2\! :=\!\omega^\omega$ and $\mathbb{A}_2\! :=\!
\big\{\big( u(2\vert u\vert )^\infty ,u(2\vert u\vert\! +\! 1)^\infty\big)\mid u\!\in\!\omega^{<\omega}\big\}\! =\!\bigcup_{i\in\omega}~C^0_i\!\times\! C^1_i$, so that $\mathbb{A}_2$ is a $\boratwo$ relation on 
$\mathbb{X}_2$. As in Section 3.(A) we get the two following consequences:
 
\begin{cor} Let $X$ be a Polish space, and $A$ be a relation on $X$. Then exactly one of the following holds:\smallskip  

\noindent (a) $(X,A)\preceq_{\bortwo}\big(\omega ,\neg\Delta (\omega )\big)$,\smallskip  

\noindent (b) $(\mathbb{X}_2,\mathbb{A}_2)\preceq_{\boraone}(X,A)$.\end{cor}
  
\begin{prop} Let $\mathbb{X}$ be a $0$-dimensional Polish space, and $\mathbb{A}$ be a relation on 
$\mathbb{X}$. The following are equivalent:\smallskip

\noindent (a) For any Polish space $X$, and any relation $A$ on $X$,
$$(X,A)\not\preceq_{\bortwo}\big(\omega ,\neg\Delta (\omega )\big)\Rightarrow 
(\mathbb{X},\mathbb{A})\preceq_{\boraone}(X,A).$$
(b) There is a comparing $1$-disjoint family $(C^\varepsilon_i)_{(\varepsilon ,i)\in 2\times\omega}$ of subsets of 
$\mathbb{X}$ such that $\mathbb{A}\!\subseteq\!\bigcup_{i\in\omega}~C^0_i\!\times\! C^1_i$.
\end{prop}

\noindent\bf (B) Countable unions of $\boratwo$ rectangles\rm

\begin{lem} (a) Assume that $(C^\varepsilon_i)_{(\varepsilon ,i)\in 2\times\omega}$ is a $1$-disjoint family of meager subsets of $\mathbb{X}$ such that no nonempty clopen subset of $\mathbb{X}$ is 
$(\bigcup_{i\in\omega}~C^0_i\!\times\! C^1_i)$-discrete. Then 
$\Delta\big(\mathbb{X}\!\setminus\! (\bigcup_{(\varepsilon ,i)\in 2\times\omega}~C^\varepsilon_i)\big)$ is not separable from $\bigcup_{i\in\omega}~C^0_i\!\times\! C^1_i$ by a 
$\Big(\boratwo\big(\mathbb{X}\!\setminus\! (\bigcup_{i\in\omega}~C^1_i)\big)\!\times\!
\boratwo\big(\mathbb{X}\!\setminus\! (\bigcup_{i\in\omega}~C^0_i)\big)\Big)_\sigma$ set.\smallskip

\noindent (b) There is a comparing $1$-disjoint family $(C^\varepsilon_i)_{(\varepsilon ,i)\in 2\times\omega}$ of subsets of $\omega^\omega$ satisfying the assumption (and thus the conclusion) of (a).\end{lem}

\noindent\bf Proof.\rm\ (a) We argue by contradiction, which gives 
$C_n\!\in\!\bormone\big(\mathbb{X}\!\setminus\! (\bigcup_{i\in\omega}~C^1_i)\big)$ on one side, and also  
$D_n\!\in\!\bormone\big(\mathbb{X}\!\setminus\! (\bigcup_{i\in\omega}~C^0_i)\big)$ with 
$\Delta\big(\mathbb{X}\!\setminus\! (\bigcup_{(\varepsilon ,i)\in 2\times\omega}~C^\varepsilon_i)\big)
\!\subseteq\!\bigcup_{n\in\omega}~(C_n\!\times\! D_n)\!\subseteq\!
\neg (\bigcup_{i\in\omega}~C^0_i\!\times\! C^1_i)$. In particular, 
$\mathbb{X}\!\setminus\! (\bigcup_{(\varepsilon ,i)\in 2\times\omega}~C^\varepsilon_i)\! =\!
\bigcup_{n\in\omega}~C_n\cap D_n$, and Baire's theorem gives $n$ and a nonempty clopen subset $C$ of $\mathbb{X}$ such that 
$C\!\setminus\! (\bigcup_{(\varepsilon ,i)\in 2\times\omega}~C^\varepsilon_i)\!\subseteq\! C_n\cap D_n$. Note that $C\!\setminus\! (\bigcup_{i\in\omega}~C^1_i)\!\subseteq\! C_n$ and 
$C\!\setminus\! (\bigcup_{i\in\omega}~C^0_i)\!\subseteq\! D_n$ since the $C^\varepsilon_i$'s are meager and $\mathbb{X}\!\setminus\! (\bigcup_{(\varepsilon ,i)\in 2\times\omega}~C^\varepsilon_i)$ is dense in $\mathbb{X}\!\setminus\! (\bigcup_{i\in\omega}~C^\varepsilon_i)$. The assumption gives 
$(x,y)\!\in\! (\bigcup_{i\in\omega}~C^0_i\!\times\! C^1_i)\cap C^2$. Then 
$(x,y)\!\in\! (\bigcup_{i\in\omega}~C^0_i\!\times\! C^1_i)\cap (C_n\!\times\! D_n)$, which is absurd.\bigskip

\noindent (b) Let $(C^\varepsilon_i)_{(\varepsilon ,i)\in 2\times\omega}$ be the family given by Lemmas 4.1.(b) and 4.3. As the $C^\varepsilon_i$'s are singletons, they are meager.\hfill{$\square$}\bigskip

\noindent\bf Remark.\rm\ Note that 
$\Delta\big(\mathbb{X}\!\setminus\! (\bigcup_{(\varepsilon ,i)\in 2\times\omega}~C^\varepsilon_i)\big)\! =\!\Delta (\mathbb{X})\cap\Big(\big(\mathbb{X}\!\setminus\! (\bigcup_{i\in\omega}~C^1_i)\big)\!\times\!
\big(\mathbb{X}\!\setminus\! (\bigcup_{i\in\omega}~C^0_i)\big)\Big)$ is a closed subset of 
$\big(\mathbb{X}\!\setminus\! (\bigcup_{i\in\omega}~C^1_i)\big)\!\times\!
\big(\mathbb{X}\!\setminus\! (\bigcup_{i\in\omega}~C^0_i)\big)$. This shows that the spaces 
$\mathbb{X}^0_2,\mathbb{X}^1_2$ of Conjecture 2 cannot be both compact, which is quite unusual in this kind of dichotomy (even if it was already the case in [L2]). Indeed, our example shows that 
$\mathbb{A}^0_2,\mathbb{A}^1_2$ must be separable by a closed set $C$, and $C,\mathbb{A}^1_2$ must have disjoint projections. If $\mathbb{X}^0_2,\mathbb{X}^1_2$ are compact, then $C$ and its projections are compact too. The product of these compact projections is a $(\boratwo\!\times\!\boratwo )_\sigma$ set separating $\mathbb{A}^0_2$ from $\mathbb{A}^1_2$, which cannot be. We will meet an example where 
$\mathbb{X}\! =\! 3^\omega$. This fact implies that we cannot extend the continuous maps of Theorem 
1.5.(e) to $3^\omega$ in general.\bigskip

 To ensure the possibility of the reduction, we introduce the following notion:
 
\begin{defi} Let $1\!\leq\!\xi\! <\!\omega_1$. A $\xi$-disjoint family 
$(C^\varepsilon_i)_{(\varepsilon ,i)\in 2\times\omega}$ of subsets of a $0$-dimensional Polish space 
$\mathbb{X}$ is said to be $very\ comparing$ if for each integer $q$ there is a partition $(O_q^p)_{p\in\omega}$ of $\mathbb{X}$ into $\borxi$ sets such that, for each $i\!\in\!\omega$,\smallskip

\noindent (a) if $q\! <\! i$, then there is $p_q^i\!\in\!\omega$ such that 
$C^0_i\cup C^1_i\!\subseteq\! O_q^{p_q^i}$,\smallskip

\noindent (b) if $q\!\geq\! i$ and $\varepsilon\!\in\! 2$, then 
$C^\varepsilon_i\!\subseteq\! O_q^{2i+\varepsilon}$,\smallskip

\noindent (c) if $(\varepsilon ,i)\!\in\! 2\!\times\!\omega$, then  
$\bigcup_{r\geq i}~\bigcap_{q\geq r}~O_q^{2i+\varepsilon}\! =\! C^\varepsilon_i$.\end{defi}

\begin{lem} There is a very comparing $1$-disjoint family 
$(C^\varepsilon_i)_{(\varepsilon ,i)\in 2\times\omega}$ of subsets of $\omega^\omega$ satisfying the assumptions (and thus the conclusion) of Lemma 4.7.(a).\end{lem}

\noindent\bf Proof.\rm\ Let $(C^\varepsilon_i)_{(\varepsilon ,i)\in 2\times\omega}$ be the family given by Lemmas 4.1.(b), 4.3, and 4.7.(b). It remains to check Condition (c). Note first that  the inclusion 
$\bigcup_{r\geq i}~\bigcap_{q\geq r}~O_q^{2i+\varepsilon}\!\supseteq\! C^\varepsilon_i$ holds for any comparing $\xi$-disjoint family.

\vfill\eject

 Conversely, 
$$\begin{array}{ll}
\bigcup_{r\geq i}~\bigcap_{q\geq r}~O_q^{2i+\varepsilon}\!\!\!\!
& \! =\!\bigcup_{r\geq i}~\bigcap_{q\geq r}~
N_{b(i)(2\vert b(i)\vert+\varepsilon )^{\mbox{max}_{l\leq q}~(\vert b(l)\vert +1)-\vert b(i)\vert}}\cr
& \! =\!\{ b(i)(2\vert b(i)\vert+\varepsilon )^\infty\}\! =\! C^\varepsilon_i.
\end{array}$$
This finishes the proof.\hfill{$\square$}\bigskip
 
\noindent\bf Notation.\rm\ We now recall some facts about the Gandy-Harrington topology (see [L2]). Let $Z$ be a recursively presented Polish space. The 
$Gandy$-$Harrington\ topology$ GH on $Z$ is generated by the $\Ana$ subsets of $Z$. We set  
$\Omega_Z :=\{ z\!\in\! Z\mid\omega_1^z\! =\!\omega_1^{\hbox{\rm CK}}\}$. Then 
$\Omega_Z$ is $\Ana$, dense in $(Z,\mbox{GH})$, and $W\cap\Omega_Z$ is a clopen subset of 
$(\Omega_Z,\mbox{GH})$ for each $W\!\in\!\Ana (Z)$. Moreover, $(\Omega_Z,\mbox{GH})$ is a $0$-dimensional Polish space. So we fix a complete compatible metric $d_{\mbox{GH}}$ on 
$(\Omega_Z,\mbox{GH})$.\bigskip

 The following notion is important for the next proof.
 
\begin{defi} Let $a$ be a countable set, $\xi\! <\!\omega_1$, and 
${\cal F}\! :=\! (S^\varepsilon_i)_{(\varepsilon ,i)\in 2\times\omega}$ be a $\xi$-disjoint family of subsets of $a^\omega$. We say that $s\!\in\! a^{<\omega}$ is ${\cal F}$-$suitable$ if there is 
$(\alpha ,\beta )\!\in\!\bigcup_{i\in\omega}~S^0_i\!\times\! S^1_i$ such that $s\! =\!\alpha\wedge\beta$ is the longest common initial segment of $\alpha$ and $\beta$.\end{defi}

\noindent\bf Example.\rm\ In the next proof, we will take $a\! :=\! 3$, $\xi\! :=\! 1$, and $s$ will be 
$suitable$ when $s$ is empty or finishes with $2$. If 
$\theta\! :\!\omega\!\rightarrow\!\{ s\!\in\! 3^{<\omega}\mid s\mbox{ is suitable}\}$ is a bijection such that 
$\big(\vert\theta (i)\vert\big)_{i\in\omega}$ is non-decreasing, then we can define a $1$-disjoint family 
$\cal F$ of subsets of $3^\omega$ by 
$S^\varepsilon_i\! :=\!\{\theta (i)\varepsilon\alpha\mid\alpha\!\in\! 2^\omega\}$, and $s$ is suitable exactly when $s$ is $\cal F$-suitable.\bigskip

 In particular, a non suitable sequence is of the form $s\varepsilon t$, where $s$ is suitable, 
$\varepsilon\!\in\! 2$ and $t\!\in\! 2^{<\omega}$ (we will use this notation in the next proof). If 
$\emptyset\!\not=\! s$ is suitable, then we set 
$$s^-\! :=\! s\vert\mbox{max}\{ l\! <\!\vert s\vert\mid s\vert l\mbox{ is suitable}\} .$$

\begin{lem} Let $\mathbb{X}$ be a $0$-dimensional Polish space, 
$(C^\varepsilon_i)_{(\varepsilon ,i)\in 2\times\omega}$ be a very comparing $1$-disjoint family of subsets of 
$\mathbb{X}$, $X,Y$ be Polish spaces, and $A_0,A_1$ be disjoint analytic subsets of $X\!\times\! Y$. Then one of the following holds:\smallskip  

\noindent (a) $A_0$ is separable from $A_1$ by a $(\boratwo\!\times\!\boratwo )_\sigma$ set,\smallskip  

\noindent (b) $\big(\mathbb{X}\!\setminus\! (\bigcup_{i\in\omega}~C^1_i),
\mathbb{X}\!\setminus\! (\bigcup_{i\in\omega}~C^0_i),
\Delta\big(\mathbb{X}\!\setminus\! (\bigcup_{(\varepsilon ,i)\in 2\times\omega}~C^\varepsilon_i)\big),
\bigcup_{i\in\omega}~C^0_i\!\times\! C^1_i\big)\leq (X,Y,A_0,A_1)$.\end{lem}

\noindent\bf Proof.\rm\ We may assume that $X,Y$ are recursively presented and that $A_0,A_1$ are 
$\Ana$. Assume that (a) does not hold. By Theorem 2.5 we get 
$N\! :=\! A_0\cap\overline{A_1}^{T_2\times T_2}\!\not=\!\emptyset$. Lemma 2.2 implies that
$$\begin{array}{ll}
(x,y)\!\notin\!\overline{A_1}^{T_2\times T_2}\!\!\!\!
& \Leftrightarrow\exists C,D\!\in\!\Ana\cap\bormone ~~(x,y)\!\in\! C\!\times\! D\!\subseteq\!\neg A_1\cr
& \Leftrightarrow\exists C,D\!\in\!\Borel\cap\bormone ~~(x,y)\!\in\! C\!\times\! D\!\subseteq\!\neg A_1.
\end{array}$$
This and Proposition 1.4 in [Lo1] show that $N$ is $\Ana$. We construct\smallskip 

\noindent - A sequence $(x_u)_{u\in 3^{<\omega}}$ of points of $X$,\smallskip

\noindent - A sequence $(y_u)_{u\in 3^{<\omega}}$ of points of $Y$,\smallskip

\noindent - A sequence $(X_u)_{u\in 3^{<\omega}}$ of $\Boraone$ subsets of $X$,\smallskip

\noindent - A sequence $(Y_u)_{u\in 3^{<\omega}}$ of $\Boraone$ subsets of $Y$,\smallskip

\noindent - A sequence $(V_s)_{s\in 3^{<\omega}\mbox{ suitable}}$ of $\Ana$ subsets of 
$X\!\times\!Y$.

\vfill\eject

 We want these objects to satisfy the following conditions:
$$\begin{array}{ll}
& (1)~(x_u,y_u)\!\in\! X_u\!\times\! Y_u\cr
& (2)~(x_s,y_s)\!\in\! V_s\!\subseteq\! N\cap\Omega_{X\times Y}\mbox{ if }s\mbox{ is suitable}\cr
& (3)~\overline{X_{u\varepsilon}}\!\subseteq\! X_u\mbox{ if }u\mbox{ is suitable or }u\! =\! s0t
\mbox{, and }\overline{X_{s1t2}}\!\subseteq\! X_s\cr
& (4)~\overline{Y_{u\varepsilon}}\!\subseteq\! Y_u\mbox{ if }u\mbox{ is suitable or }u\! =\! s1t
\mbox{, and }\overline{Y_{s0t2}}\!\subseteq\! Y_s\cr
& (5)~V_s\!\subseteq\! V_{s^-}\mbox{ if }\emptyset\!\not=\! s\mbox{ is suitable}\cr
& (6)~\mbox{diam}(X_u),\mbox{diam}(Y_u)\!\leq\! 2^{-\vert u\vert}\cr
& (7)~\mbox{diam}_{\mbox{GH}}(V_s)\!\leq\! 2^{-\vert s\vert}\mbox{ if }s\mbox{ is suitable}\cr
& (8)~(x_{s0},y_{s1})\!\in\!\big(\overline{\Pi_0[(X_s\!\times\! Y_s)\cap V_s]}\!\times\!
\overline{\Pi_1[(X_s\!\times\! Y_s)\cap V_s]}\big)\cap A_1\mbox{ if }s\mbox{ is suitable}\cr
& (9)~(x_{s0t},y_{s1t})\! =\! (x_{s0},y_{s1})\mbox{ if }s\mbox{ is suitable and }t\!\in\! 2^{<\omega}
\end{array}$$ 
\noindent $\bullet$ Assume that this is done. We define 
$\phi\! :\!\omega^{<\omega}\!\rightarrow\! 3^{<\omega}$ by $\phi (\emptyset )\! :=\!\emptyset$ and
$$\phi (sn)\! :=\!\left\{\!\!\!\!\!\!\!
\begin{array}{ll}
& \phi (s)2\varepsilon\mbox{ if }n\! =\! 2\vert s\vert\! +\!\varepsilon\mbox{ or }
\big( s\!\not=\!\emptyset\mbox{ and }n\!\not=\! s(\vert s\vert\! -\! 1)\big)\mbox{,}\cr
& \phi (s)\varepsilon\mbox{ if }(n\! =\! 2q\! +\!\varepsilon\mbox{ and }q\!\not=\!\vert s\vert)\mbox{ and }
\Big( s\! =\!\emptyset\mbox{ or }\big( s\!\not=\!\emptyset\mbox{ and }n\! =\! s(\vert s\vert\! -\! 1)\big)\Big).
\end{array}
\right.$$
This map allows us to define $\Phi\! :\!\omega^{\omega}\!\rightarrow\! 3^{\omega}$ by 
$\Phi (\gamma )(p)\! :=\!\phi\big(\gamma\vert (p\! +\! 1)\big)(p)$, and $\Phi$ is continuous.\bigskip

 As $(C^\varepsilon_i)_{(\varepsilon ,i)\in 2\times\omega}$ is very comparing, there are some witnesses 
$(O_q^p)_{p\in\omega}$. Let $x\!\in\!\mathbb{X}$. As in the proof of Lemma 4.4, we associate the sequence $(p_q)_{q\in\omega}\!\in\!\omega^\omega$ defined by $x\!\in\! O^{p_q}_q$. As 
$(C^\varepsilon_i)_{(\varepsilon ,i)\in 2\times\omega}$ is very comparing, $(p_q)_{q\in\omega}$ is not eventually constant if $x\!\notin\!\bigcup_{(\varepsilon ,i)\in 2\times\omega}~C^\varepsilon_i$. Thus 
$\Phi\big( (p_q)_{q\in\omega}\big)$ has infinitely many $2$'s in this case. If $x\!\in\! C^\varepsilon_i$, then 
$$\Phi\big( (p_q)_{q\in\omega}\big)\! =\!\Phi\big( p^i_0...p^i_{i-1}(2i\! +\!\varepsilon )^\infty\big)\! =\!
\phi (p^i_0...p^i_{i-1})2\varepsilon^\infty .$$
If $x\!\in\!\mathbb{X}\!\setminus\! (\bigcup_{i\in\omega}~C^1_i)$, then the increasing sequence 
$(n^0_k)_{k\in\omega}$ of integers such that $\Phi\big( (p_q)_{q\in\omega}\big)\vert n^0_k$ is suitable or of the form $s0t$ is infinite. Condition (3) implies that 
$(\overline{X_{\Phi ((p_q)_{q\in\omega})\vert n^0_k}})_{k\in\omega}$ is non-increasing. Moreover, 
$(\overline{X_{\Phi ((p_q)_{q\in\omega})\vert n^0_k}})_{k\in\omega}$ is a sequence of nonempty closed subsets of $X$ whose diameters tend to $0$, so that we can define 
$\{ f(x)\}\! :=\!\bigcap_{k\in\omega}~\overline{X_{\Phi ((p_q)_{q\in\omega})\vert n^0_k}}\! =\!
\bigcap_{k\in\omega}~X_{\Phi ((p_q)_{q\in\omega})\vert n^0_k}$. This defines a continuous map 
$f\! :\!\mathbb{X}\!\setminus\! (\bigcup_{i\in\omega}~C^1_i)\!\rightarrow\! X$ with 
$f(x)\! =\!\mbox{lim}_{k\rightarrow\infty}~x_{\Phi ((p_q)_{q\in\omega})\vert n^0_k}$. Similarly, we define 
$g\! :\!\mathbb{X}\!\setminus\! (\bigcup_{i\in\omega}~C^0_i)\!\rightarrow\! Y$ continuous with 
$g(x)\! =\!\mbox{lim}_{k\rightarrow\infty}~y_{\Phi ((p_q)_{q\in\omega})\vert n^1_k}$.\bigskip

 If $x\!\notin\!\bigcup_{(\varepsilon ,i)\in 2\times\omega}~C^\varepsilon_i$, then the sequence 
$(k_j)_j$ of integers such that $\Phi\big( (p_q)_{q\in\omega}\big)\vert k_j$ is suitable is an infinite subsequence of both $(n^0_k)_{k\in\omega}$ and $(n^1_k)_{k\in\omega}$. Note that 
$(V_{\Phi ((p_q)_{q\in\omega})\vert k_j})_{j\in\omega}$ is a non-increasing sequence of nonempty closed subsets of $\Omega_{X\times Y}$ whose GH-diameters tend to $0$, so that we can define $F(x)$ by 
$\{ F(x)\}\! :=\!\bigcap_{j\in\omega}~V_{\Phi ((p_q)_{q\in\omega})\vert k_j}\!\subseteq\! N\!\subseteq\! A_0$. As $F(x)$ is the limit (in $(X\!\times\! Y,\mbox{GH})$, and thus in $X\!\times\! Y$) of 
$(x_{\Phi ((p_q)_{q\in\omega})\vert k_j},y_{\Phi ((p_q)_{q\in\omega})\vert k_j})_{j\in\omega}$, we get 
$F(x)\! =\!\big( f(x),g(x)\big)$. Therefore  
$\Delta\big(\mathbb{X}\!\setminus\! (\bigcup_{(\varepsilon ,i)\in 2\times\omega}~C^\varepsilon_i)\big)
\!\subseteq\! (f\!\times\!g)^{-1}(A_0)$.\bigskip

 Note that $x_{s0}\! =\! x_{s0^2}\! =\! ...\! =\! x_{s0^{q+1}}$ for each $s$ suitable. Thus 
$$f(x)\! =\!\mbox{lim}_{q\rightarrow\infty}~x_{\phi (p^i_0...p^i_{i-1})20^{q+1}}\! =\! 
x_{\phi (p^i_0...p^i_{i-1})20}$$ 
if $x\!\in\! C^0_i$. Similarly, $g(y)\! =\! y_{\phi (p^i_0...p^i_{i-1})21}$ if $y\!\in\! C^1_i$ and 
$\bigcup_{i\in\omega}~C^0_i\!\times\! C^1_i\!\subseteq\! (f\!\times\!g)^{-1}(A_1)$.

\vfill\eject

\noindent $\bullet$ Let us prove that the construction is possible. As $N\!\not=\!\emptyset$, we can choose 
$(x_\emptyset ,y_\emptyset )\!\in\! N\cap\Omega_{X\times Y}$, a $\Ana$ subset $V_\emptyset$ of 
$X\!\times\! Y$ with $(x_\emptyset ,y_\emptyset )\!\in\! V_\emptyset\!\subseteq\! N\cap\Omega_{X\times Y}$ of GH-diameter at most $1$, and a $\Boraone$ neighborhood $X_\emptyset$ (resp., $Y_\emptyset$) of 
$x_\emptyset$ (resp., $y_\emptyset$) of diameter at most $1$. Assume that $(x_u)_{u\in 3^{\leq l}}$, 
$(y_u)_{u\in 3^{\leq l}}$, $(X_u)_{u\in 3^{\leq l}}$, $(Y_u)_{u\in 3^{\leq l}}$ and 
$(V_s)_{s\in 3^{\leq l}\mbox{ suitable}}$ satisfying (1)-(9) have been constructed, which is the case for 
$l\! =\! 0$.\bigskip

 Let $s\!\in\! 3^{<\omega}$ be suitable. Note that 
$(x_s,y_s)\!\in\! (X_s\!\times\! Y_s)\cap V_s\!\subseteq\!\overline{A_1}^{T_2\times T_2}$. We choose 
$X',Y'\!\in\!\Boraone$ with 
$(x_s,y_s)\!\in\! X'\!\times\! Y'\!\subseteq\!\overline{X'}\!\times\!\overline{Y'}\!\subseteq\! X_s\!\times\! Y_s$. As $\Pi_\varepsilon [(X'\!\times\! Y')\cap V_s]$ is $\Ana$, 
$\overline{\Pi_\varepsilon [(X'\!\times\! Y')\cap V_s]}$ is $\Ana\cap\bormone$. In particular, 
$\overline{\Pi_\varepsilon [(X'\!\times\! Y')\cap V_s]}$ is $T_2$-open. This shows the existence of 
$$(x_{s0},y_{s1})\!\in\!\big(\overline{\Pi_0[(X'\!\times\! Y')\cap V_s]}\!\times\!
\overline{\Pi_1[(X'\!\times\! Y')\cap V_s]}\big)\cap A_1.$$ 
Note that $(x_{s0},y_{s1})\!\in\!\overline{X'}\!\times\!\overline{Y'}\!\subseteq\! X_s\!\times\! Y_s$. We set 
$x_{s1}\! :=\! x_s$, $y_{s0}\! :=\! y_s$. We defined $x_u,y_u$ when $u\!\in\! 3^{l+1}$ is not suitable but 
$u\vert l$ is suitable.\bigskip

 Assume now that $u\!\in\! 3^{l+1}$ is suitable, but not $u\vert l$. This gives $(s,\varepsilon ,t)$ such that 
$u\! =\! s\varepsilon t2$. Assume first that $\varepsilon\! =\! 0$. Note that 
$x_{s0t}\! =\! x_{s0}\!\in\! X_{s0t}\cap\overline{\Pi_0[(X_s\!\times\! Y_s)\cap V_s]}$. This gives 
$$x_u\!\in\! X_{s0t}\cap\Pi_0[(X_s\!\times\! Y_s)\cap V_s]\mbox{,}$$ 
and also $y_u$ with $(x_u,y_u)\!\in\!\big( (X_{s0t}\cap X_s)\!\times\! Y_s\big)\cap V_s\! =\! 
(X_{s0t}\!\times\! Y_s)\cap V_s$. If $\varepsilon\! =\! 1$, then similarly we get 
$(x_u,y_u)\!\in\! (X_s\!\times\! Y_{s1t})\cap V_s$.\bigskip

 If $u$ and $u\vert l$ are both suitable, or both non suitable, then we set 
$(x_u,y_u)\! :=\! (x_{u\vert l},y_{u\vert l})$. So we defined $x_u,y_u$ in any case. Note that Conditions 
(8) and (9) are fulfilled, and that $(x_s,y_s)\!\in\! V_{s^-}$ if $s$ is suitable. Moreover, 
$x_u\!\in\! X_{u\vert l}$ if $u\vert l$ is suitable or $u\vert l\! =\! s0t$, and $x_u\!\in\! X_s$ if 
$u\! =\! s1t2$, and similarly in $Y$. We choose $\Boraone$ sets $X_u,Y_u$ of diameter at most 
$2^{-l-1}$ with
$$(x_u,y_u)\!\in\! X_u\!\times\! Y_u\!\subseteq\!\overline{X_u}\!\times\!\overline{Y_u}\!\subseteq\!
\left\{\!\!\!\!\!\!
\begin{array}{ll}
& X_{u\vert l}\!\times\! Y_{u\vert l}\mbox{ if }u\mbox{ is not suitable or }u\vert l\mbox{ is suitable,}\cr
& X_{u\vert l}\!\times\! Y_s\mbox{ if }u\! =\! s0t2\mbox{,}\cr
& X_s\!\times\! Y_{u\vert l}\mbox{ if }u\! =\! s1t2.
\end{array}
\right.$$
It remains to choose, when $s$ is suitable, $V_s\!\in\!\Ana (X\!\times\! Y)$ of GH-diameter at most 
$2^{-l-1}$ such that $(x_s,y_s)\!\in\! V_s\!\subseteq\! V_{s^-}$.\hfill{$\square$}\bigskip

 We set $\mathbb{X}^\varepsilon_2\! :=\!
\omega^\omega\!\setminus\!\{ u(2\vert u\vert\! +\! 1\! -\!\varepsilon )^\infty\mid u\!\in\!\omega^{<\omega}\}\! =\!
\mathbb{X}\!\setminus\! (\bigcup_{i\in\omega}~C^{1-\varepsilon}_i)$, 
$$\mathbb{A}^0_2\! :=\!\Delta (\omega^\omega\!\setminus\!
\{ u(2\vert u\vert\! +\!\varepsilon )^\infty\mid (u,\varepsilon )\!\in\!\omega^{<\omega}\!\times\! 2\})\! =\!
\Delta\big(\mathbb{X}\!\setminus\! (\bigcup_{(\varepsilon ,i)\in 2\times\omega}~C^\varepsilon_i)\big)
\mbox{,}$$
and $\mathbb{A}^1_2\! :=\!
\big\{\big( u(2\vert u\vert )^\infty ,u(2\vert u\vert\! +\! 1)^\infty\big)\mid u\!\in\!\omega^{<\omega}\big\}\! =\!\bigcup_{i\in\omega}~C^0_i\!\times\! C^1_i$.

\begin{cor} Let $X,Y$ be Polish spaces, and $A_0,A_1$ be disjoint analytic subsets of $X\!\times\! Y$. Then exactly one of the following holds:\smallskip  

\noindent (a) $A_0$ is separable from $A_1$ by a $(\boratwo\!\times\!\boratwo )_\sigma$ set,\smallskip  

\noindent (b) $(\mathbb{X}^0_2,\mathbb{X}^1_2,\mathbb{A}^0_2,\mathbb{A}^1_2)\leq (X,Y,A_0,A_1)$.
\end{cor}

\vfill\eject

\noindent\bf Remark.\rm\ In the remark after Lemma 4.7 we announced an example with 
$\mathbb{X}\! :=\! 3^\omega$. In fact, we already met it after Definition 4.10. Recall that the formula 
$S^\varepsilon_i\! :=\!\{\theta (i)\varepsilon\alpha\mid\alpha\!\in\! 2^\omega\}$ defines a $1$-disjoint family of subsets of $3^\omega$, which are clearly meager. It is also clear that no nonempty clopen subset of $3^\omega$ is $(\bigcup_{i\in\omega}~S^0_i\!\times\! S^1_i)$-discrete. One can check that the formula
$$O^p_q\! :=\!\left\{\!\!\!\!\!\!\!
\begin{array}{ll}
& \bigcup_{t\in 2^{\vert\theta (q)\vert -\vert\theta (i)\vert}}~N_{\theta (i)\varepsilon t}
\mbox{ if }p\! =\! 2i\! +\!\varepsilon\!\leq\! 2q\! +\! 1\mbox{,}\cr\cr
& 3^\omega\!\setminus\! (\bigcup_{p'\leq 2q+1}~O^{p'}_q)\mbox{ if }p\! =\! 2q\! +\! 2\mbox{,}\cr\cr
& \emptyset\mbox{ if }p\!\geq\! 2q\! +\! 3\mbox{,}
\end{array}
\right.$$
defines witnesses for the fact that $(S^\varepsilon_i)_{(\varepsilon ,i)\in 2\times\omega}$ is very comparing.\bigskip

 To close this section, we notice that the notion of a very comparing $1$-disjoint family gives only a sufficient condition, and not a characterization like in 3.5, 3.9 or 4.6:
   
\begin{prop} Let $\mathbb{X}$ be a $0$-dimensional Polish space, $(C^\varepsilon_i)_{i\in\omega}$ be a very comparing $1$-disjoint family of subsets of $\mathbb{X}$, $\mathbb{X}^\varepsilon\!\subseteq\!\mathbb{X}\!\setminus\! (\bigcup_{i\in\omega}~C^{1-\varepsilon}_i)$, 
$\mathbb{A}^0\!\subseteq\!
\Delta\big(\mathbb{X}\!\setminus\! (\bigcup_{(\varepsilon ,i)\in 2\times\omega}~C^\varepsilon_i)\big)$, and $\mathbb{A}^1\!\subseteq\!\bigcup_{i\in\omega}~C^0_i\!\times\! C^1_i$ be as in the definition of 
$\leq$. Then for any Polish spaces $X,Y$, and any disjoint analytic subsets $A_0,A_1$ of 
$X\!\times\! Y$,
$$A_0\mbox{ is not separable from }A_1\mbox{ by a }(\boraone\!\times\!\boraone )_\sigma\mbox{ set }\Rightarrow 
(\mathbb{X}^0,\mathbb{X}^1,\mathbb{A}^0,\mathbb{A}^1)\leq (X,Y,A_0,A_1).$$
\end{prop}

\section{$\!\!\!\!\!\!$ The case $\xi\! =\! 3$: Baire class two colorings}

\bf Remark.\rm\ Unlike when $\xi\!\in\!\{ 1,2\}$, we cannot have $\mathbb{A}_3$ of the form 
$\bigcup_{n\in\omega}~C^0_n\!\times\! C^1_n$, where 
$(C^\varepsilon_n)_{(\varepsilon ,n)\in 2\times\omega}$ is a $2$-disjoint family. Indeed, we will see that there is a Borel graph $\mathbb{G}\!\subseteq\! 2^\omega\!\times\! 2^\omega$ of a partial injection such that $(2^\omega ,\mathbb{G})\not\preceq_{\borthree}\big(\omega ,\neg\Delta (\omega )\big)$. We would 
get $f\! :\!\mathbb{X}_3\!\rightarrow\! 2^\omega$ continuous such that 
$\mathbb{A}_3\!\subseteq\! (f\!\times\! f)^{-1}(\mathbb{G})$, and $f[C^0_n]\!\times\! f[C^1_n]$ would be a singleton. The set $(f\!\times\! f)[\mathbb{A}_3]$ would be countable, and  
$(2^\omega ,(f\!\times\! f)[\mathbb{A}_3])\preceq_{\borthree}\big(\omega ,\neg\Delta (\omega )\big)$, 
$(\mathbb{X}_3,\mathbb{A}_3)\preceq_{\borthree}\big(\omega ,\neg\Delta (\omega )\big)$ would hold, which is absurd. However, the following result holds.

\begin{thm} There are a $0$-dimensional Polish space $\mathbb{X}_3$ and an analytic relation 
$\mathbb{A}_3$ on $\mathbb{X}_3$ such that for any Polish space $X$, and for any analytic relation $A$ on $X$, exactly one of the following holds:\smallskip  

\noindent (a) $(X,A)\preceq_{\borthree}\big(\omega ,\neg\Delta (\omega )\big)$,\smallskip  

\noindent (b) $(\mathbb{X}_3,\mathbb{A}_3)\preceq_{\boraone}(X,A)$.\end{thm}

 We can take $\mathbb{X}_3\! =\!\omega^\omega$, but this is not the most natural thing to do. Note that we can replace $\mathbb{X}_3$ with any copy of it. Our space $\mathbb{X}_3$ will be a dense $G_\delta$ subset of $2^\omega$, in fact a copy of $\omega^\omega$. This $G_\delta$ subset is not necessary to see that $(\mathbb{X}_3,\mathbb{A}_3)$ satisfies the ``exactly" part of Theorem 5.1 (i.e., that $(\mathbb{X}_3,\mathbb{A}_3)\not\preceq_{\borthree}\big(\omega ,\neg\Delta (\omega )\big)$), but it is useful to build and ensure the continuity of the homomorphism of Statement (b). The definition of 
 $\mathbb{X}_3$ and $\mathbb{A}_3$ is based on the construction of the following basic objects.\bigskip
 
 Recall the sequence $(s_n)_{n\in\omega}$ defined in the introduction. In $\mathbb{G}_0$, we put $s_n$, i.e., a finite sequence of elements of $2$, before the changed coordinate. In $\mathbb{A}_3$, we will put a finite sequence of elements of $2^{<\omega}$, together with a way to recover them after concatenation, before the changed coordinate. In order to do that, we identify $\omega$ with $\omega^2$.\bigskip  
 
\noindent\bf Notation.\rm\ Let $<.,.>:\!\omega^2\!\rightarrow\!\omega$ be a natural bijection. 
More precisely, $<\! n,p\! >:=\! (\Sigma_{k\leq n+p}~k)\! +\! p$. Note that the inverse bijection 
$q\!\mapsto\!\big((q)_{0},(q)_{1}\big)$ is build as follows. We set, for $q\!\in\!\omega$,
$$M(q)\! :=\!\hbox{\rm max}\{ m\!\in\!\omega\mid\Sigma_{k\leq m}~k\!\leq\! q\}.$$
Then we define $\big((q)_{0},(q)_{1}\big)\! :=\!\big(M(q)\! -\! q\! +\! 
(\Sigma_{k\leq M(q)}~k),q\! -\! (\Sigma_{k\leq M(q)}~k)\big)$. More concretely,
$$\omega\! =\!\{<0,0>,<1,0>,<0,1>,\ldots ,<M(q),0>,<M(q)\! -\! 1,1>,...,<0,M(q)>,...\}.$$
If $u\!\in\! 2^{\leq\omega}$ and $n\!\in\!\omega$, then we define $(u)_n\!\in\! 2^{\leq\omega}$ by 
$(u)_n(p)\! :=\! u(<n,p>)$ if $<n,p><\!\vert u\vert$. Here also we define 
$<\alpha_0,\alpha_1,...>\in\! 2^\omega$ by $<\alpha_0,\alpha_1,...>(<n,p>)\! :=\!\alpha_n(p)$, for any sequence $(\alpha_n)_{n\in\omega}$ of elements of $2^\omega$. In particular, 
$\alpha\!\mapsto\!\big( (\alpha )_n\big)_{n\in\omega}$ and 
$(\alpha_n )_{n\in\omega}\!\mapsto <\alpha_0,\alpha_1,...>$ are inverse bijections.

\begin{lem} Let $u,v\!\in\! 2^{<\omega}$.\smallskip

\noindent (a) $u\!\subseteq\! v$ implies that $(u)_n\!\subseteq\! (v)_n$ for each $n\!\in\!\omega$.\smallskip

\noindent (b) $\vert (u)_0\vert\!\leq\!\vert u\vert$.\smallskip

\noindent (c) $\vert (u)_n\vert\!\leq\!\vert u\vert\! +\! 1\! -\! n$ if $n\!\leq\!\vert u\vert\! +\! 1$.\end{lem}

\noindent\bf Proof.\rm\ (a) If $<n,p><\!\vert u\vert$, then $(u)_n(p)\! =\! u(<n,p>)\! =\! v(<n,p>)\! =\! (v)_n(p)$ because of the inequality $<n,p><\!\vert v\vert$, so that $(u)_n\!\subseteq\! (v)_n$.\bigskip

\noindent (b) We set, for $n,q\!\in\!\omega$, $c^n_q\! :=\!\mbox{Card}(\{ p\!\in\!\omega\mid<n,p><\! q\})$. As $<.,.>$ is a bijection, we get $c^n_{q+1}\!\leq\! c^n_q\! +\! 1$. As $c^n_0\! =\! 0$, $c^n_q\!\leq\! q$. We are done since $\vert (u)_n\vert\! =\! c^n_{\vert u\vert}$.\bigskip

\noindent (c) Note first that 
$<\! n,p\! >=\! (\Sigma_{k\leq n+p}~k)\! +\! p\! <\! (\Sigma_{k\leq n'+p'}~k)\! +\! p'\! =<n',p'>$ if 
$n\! +\! p\! <\! n'\! +\! p'$, and that $(q)_{0}\! +\! (q)_{1}\! =\! M(q)\!\leq\! q\! <\! q\! +\! 1$. This implies that 
$q\! =<(q)_{0},(q)_{1}><<n,q\! +\! 1\! -\! n>$ if $n\!\leq\! q\! +\! 1$. It remains to apply this to 
$q\! :=\!\vert u\vert$ since $\vert (u)_n\vert\! =\! c^n_{\vert u\vert}$.\hfill{$\square$}\bigskip

 We can view $\mathbb{G}_0$ as the countable union $\bigcup_{n\in\omega}~\mbox{Gr}(\varphi_n)$, where $\varphi_n$ is the homeomorphism defined on the basic clopen set $N_{s_n0}$ onto the clopen set 
$N_{s_n1}$ defined by $\varphi_n(s_n0\gamma )\! :=\! s_n1\gamma$. The set $\mathbb{A}_3$ will also be the countable union of the graphs of some homeomorphisms, indexed by 
$\omega^{<\omega}$ instead of $\omega$. Their domain and range will be $G_\delta$ subsets of 
$2^\omega$ instead of clopen sets. We first define the closures of these $G_\delta$'s. They will be copies of $2^\omega$. In fact, our homeomorphims will also be defined on the closure of these final domains. We will fix the coordinates whose number is in one of the verticals before that of the number of the changed coordinate. This leads to the following notation.\bigskip
 
\noindent\bf Notation.\rm\ If $t\!\in\!\omega^{<\omega}$ and $k\!\leq\!\vert t\vert$, then we set 
$\Sigma^t_k\! :=\!\Sigma_{j<k}~\big( t(j)\! +\! 2\big)$, and 
$$\Sigma_t\! :=<\Sigma_{j<\vert t\vert}~\big( t(j)\! +\! 2\big),0>=<\Sigma^t_{\vert t\vert},0>$$
($\Sigma_t$ will be the number of the unique changed coordinate). We set $w_n\! :=\! s_n0$, so that 
$\vert w_n\vert\! =\! n\! +\! 1$ and $(w_n)_{n\in\omega}$ is dense (we want $w_n$ to be nonempty).

\vfill\eject

 We define the following objects for $t\!\in\!\omega^{<\omega}$.\bigskip

\noindent $\bullet$ We first define a copy $K_t$ of $2^\omega$ by\bigskip

\leftline{$K_t\! :=\!\big\{\alpha\!\in\! 2^\omega\mid\forall k\! <\!\vert t\vert ~~(\alpha )_{\Sigma^t_k}\! =\!
(w_{t(k)})_00^{t(k)+1-\vert (w_{t(k)})_0\vert}10^\infty\mbox{ and }$}\smallskip

\rightline{$\forall 0\! <\! i\! <\! t(k)\! +\! 2~~(\alpha )_{\Sigma^t_k+i}\! =\!\big( w_{t(k)}\big)_i0^\infty\big\}$.}\bigskip

\noindent This is well defined since $\vert w_{t(k)}\vert\! =\! t(k)\! +\! 1$, so that we can apply Lemma 5.2.(b) to $u\! :=\! w_{t(k)}$ and $t(k)\! +\! 1\! -\!\vert (w_{t(k)})_0\vert\!\geq\! 0$. In particular, the last $1$ in 
$(\alpha )_{\Sigma^t_k}$ is at the position $t(k)\! +\! 1$. Here is the picture of $K_t$ when $t\! =\! (4,2)$:\bigskip

\includegraphics[scale=0.552]{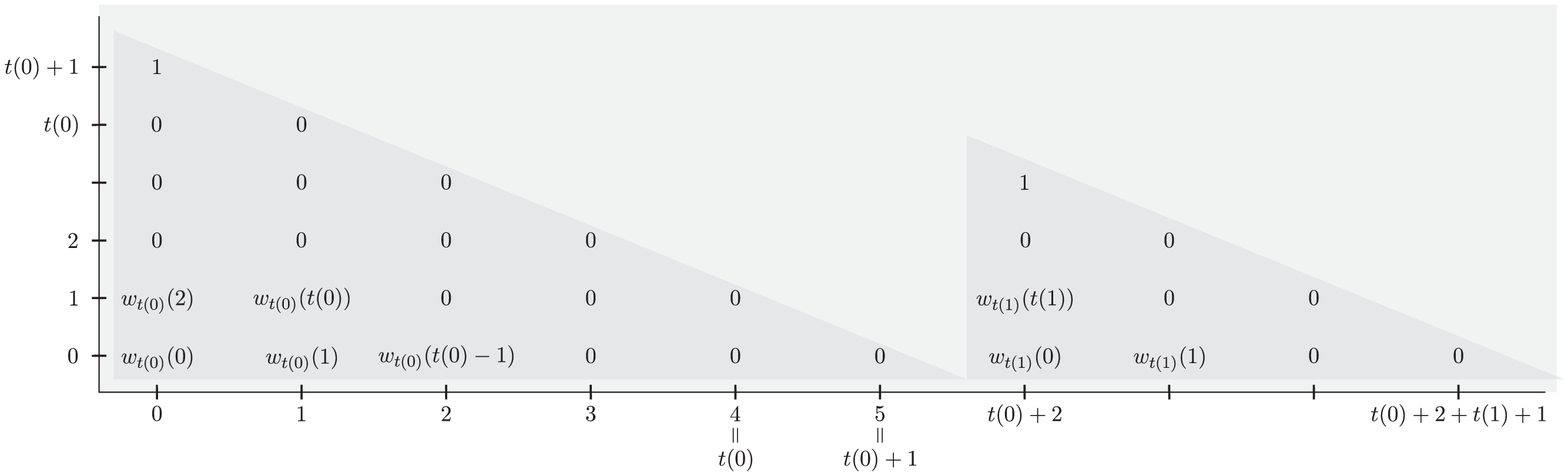}\smallskip

\noindent $\bullet$ We define a non-trivial partition $(K^0_t,K^1_t)$ of $K_t$ into clopen sets by
$K^\varepsilon_t\! :=\!\{\alpha\!\in\! K_t\mid\alpha (\Sigma_t)\! =\!\varepsilon\}$.\bigskip
 
\noindent $\bullet$ We define a homeomorphism $\varphi_t\! :\! K^0_t\!\rightarrow\! K^1_t$ by 
$\varphi_t(\alpha )(m)\! :=\! 1$ if $m\! =\!\Sigma_t$, $\alpha (m)$ otherwise.\bigskip

 We can view the construction of $K_t$, $K^\varepsilon_t$ and $\varphi_t$ inductively. Indeed, 
$K_\emptyset\! =\! 2^\omega$, $K^\varepsilon_\emptyset$ is the basic clopen set $N_\varepsilon$, and 
$\varphi_\emptyset (\alpha )(m)$ is $1$ when $m\! =\! 0$, $\alpha (m)$ otherwise. Then
$$K_{tn}\! :=\!\{\alpha\!\in\! K^{w_n(0)}_t\mid
(\alpha )_{\Sigma^t_{\vert t\vert}}\! =\! (w_n)_00^{n+1-\vert (w_n)_0\vert}10^\infty
\mbox{ and }\forall 0\! <\! i\! <\! n\! +\! 2~~(\alpha )_{\Sigma^t_{\vert t\vert}+i}\! =\! (w_n)_i0^\infty\}\mbox{,}$$
$K^\varepsilon_{tn}\! :=\!
\{\alpha\!\in\! K_{tn}\mid\alpha (<\Sigma^t_{\vert t\vert}\! +\! n\! +\! 2,0>)\! =\!\varepsilon\}$, and 
$\varphi_{tn}(\alpha )(m)$ is $1$ when $m$ is equal to $<\Sigma^t_{\vert t\vert}\! +\! n\! +\! 2,0>$, 
$\alpha (m)$ otherwise.\bigskip

 The set $\mathbb{S}_3\! :=\!\{\alpha\!\in\! 2^\omega\mid\exists m\!\in\!\omega ~~\forall n\!\in\!\omega ~~\exists p\!\geq\! n~~(p)_0\! =\! m\mbox{ and }\alpha (p)\! =\! 1\}$ is a standard $\borathree$-complete set (see 23.A in [K]). We will more or less recover this example, but the $1$'s have to be well placed. This leads to the following technical but crucial notion.

\begin{defi} We say that $u\!\in\! 2^{<\omega}$ is $placed$ if $u\!\not=\!\emptyset$ and there is 
$t\!\in\!\omega^{<\omega}$ such that $N_u\cap K_t\!\not=\!\emptyset$, 
$(\vert u\vert\! -\! 1)_0\! =\!\Sigma^t_{\vert t\vert}$, and $u(\vert u\vert\! -\! 1)\! =\! 1$ if 
$(\vert u\vert\! -\! 1)_1\! >\! 0$. We also say  that $t$ is a $witness$ for the fact that $u$ is placed.\end{defi}

 This means that the last coordinate of $u$ has a number on the vertical $\Sigma^t_{\vert t\vert}$, on which the coordinates of the elements of $K_t$ are left free by $t$, and which is the first vertical with this property. The coordinates of $u$ whose number is on one of the verticals before the previous one are determined by $t$. Finally, the last coordinate of $u$ is $1$, except maybe if this coordinate has the number $\Sigma_t$, which is at the bottom of the vertical $\Sigma^t_{\vert t\vert}$.
 
\vfill\eject
 
\noindent\bf Examples.\rm\ Let $\alpha\!\in\! K_{nj}\! =\! K_{(n,j)}$. Then\smallskip

\noindent $\bullet$ $\alpha\vert 1$, $\alpha\vert (<\! 0,n\! +\! 1\! >\! +1)$ are placed with witness 
$\emptyset$.\smallskip

\noindent $\bullet$ $\alpha\vert (<\! n\! +\! 2,0\! >\! +1)$, $\alpha\vert (<\! n\! +\! 2,j\! +\! 1\! >\! +1)$ are placed with witness $(n)$.\smallskip

\noindent $\bullet$ $\alpha\vert (<\! n\! +\! 2\! +\! j\! +\! 2,0\! >\! +1)$ is placed with witness $(n,j)$. If 
$\alpha (<\! n\! +\! 2\! +\! j\! +\! 2,q\! >)\! =\! 1$, then $\alpha\vert (<\! n\! +\! 2\! +\! j\! +\! 2,q\! >\! +1)$ is placed with witness $(n,j)$.\bigskip

 We are now ready to define $\mathbb{X}_3$ and $\mathbb{A}_3$.\bigskip

\noindent\bf Notation.\rm\ We set $\mathbb{X}_3\! :=\!\big\{\alpha\!\in\! 2^\omega\mid\forall n\!\in\!\omega ~~
\exists p\!\geq\! n~~\alpha\vert p\mbox{ is placed}\big\}$. Let $t\!\in\!\omega^{<\omega}$. We set 
$$H_t\! :=\!\{\alpha\!\in\! K^0_t\mid\forall n\!\in\!\omega ~~\exists p\!\geq\! n~~
(p)_0\! =\!\Sigma^t_{\vert t\vert}\mbox{ and }\alpha (p)\! =\! 1\}\mbox{,}$$
and $\mathbb{A}_3\! :=\!\bigcup_{t\in\omega^{<\omega}}~\mbox{Gr}({\varphi_t}_{\vert H_t})$. In this sense, we recover $\mathbb{S}_3$. More concretely,
$$\mathbb{A}_3\! =\!\bigcup_{t\in\omega^{<\omega}}~\Big\{ (u0\gamma ,u1\gamma )\mid
\vert u\vert\! =\!\Sigma_t\mbox{ and }u0\gamma\!\in\! K_t\mbox{ and }\forall n\!\in\!\omega ~~
\exists p\!\geq\! n~~(u0\gamma )(<\Sigma^t_{\vert t\vert},p>)\! =\! 1\Big\} .$$
 
\begin{lem} Let $t\!\in\!\omega^{<\omega}$ and $\varepsilon\!\in\! 2$.\smallskip

\noindent (a) $(K_{tn})_{n\in\omega ,w_n(0)=\varepsilon}$ is a sequence of pairwise disjoint meager subsets of $K^\varepsilon_t$.\smallskip

\noindent (b) Any nonempty open subset of $K^\varepsilon_t$ contains one of the $K_{tn}$'s.\end{lem}

\noindent\bf Proof.\rm\ (a) This comes from the fact that the last $1$ in 
$(\alpha )_{\Sigma^t_{\vert t\vert}}$ is at the position $n\! +\! 1$ if $\alpha\!\in\! K_{tn}$.\bigskip

\noindent (b) A nonempty open subset of $K^\varepsilon_t$ contains a basic clopen set $C$ of the form 
$$\{\alpha\!\in\! K^\varepsilon_t\mid\varepsilon u\!\subseteq <(\alpha )_{\Sigma^t_{\vert t\vert}},
(\alpha )_{\Sigma^t_{\vert t\vert}+1},...>\}\mbox{,}$$
where $u\!\in\! 2^{<\omega}$. We choose $n\!\in\!\omega$ such that $\varepsilon u\!\subseteq\! w_n$. It remains to see that $K_{tn}\!\subseteq\! C$. So let 
$m\! =<i,p>\leq\vert u\vert$. Note first that $M(q)\!\leq\!\mbox{min}\big( q,M(q\! +\! 1)\big)$. Thus, as  
$\vert u\vert\!\leq\!\vert w_n\vert\! =\! n\! +\! 1$, 
$$i\! =\! (m)_0\!\leq\! (m)_0\! +\! (m)_1\! =\! M(m)\!\leq\! 
M(\vert u\vert )\!\leq\! M(n\! +\! 1)\!\leq\! n\! +\! 1\! <\! n\! +\! 2.$$
Lemma 5.2.(a) allows us to write
$$\begin{array}{ll}
(\varepsilon u)(m)\!\!\!
& \! =\! (\varepsilon u)(<i,p>)\! =\! (\varepsilon u)_i(p)\! =\! (w_n)_i(p)\! =\! 
(\alpha )_{\Sigma^t_{\vert t\vert}+i}(p)\cr
& \! =\! (<(\alpha )_{\Sigma^t_{\vert t\vert}},(\alpha )_{\Sigma^t_{\vert t\vert}+1},...>)_i(p)
\! =\! (<(\alpha )_{\Sigma^t_{\vert t\vert}},(\alpha )_{\Sigma^t_{\vert t\vert}+1},...>)(m).
\end{array}$$
This finishes the proof.\hfill{$\square$}\bigskip

 We now start to prove the required properties of $\mathbb{X}_3$ and $\mathbb{A}_3$.

\begin{lem} (a) The set $\mathbb{X}_3$ is a dense $\bormtwo$ subset of $2^\omega$. In particular, 
$\mathbb{X}_3$ is a $0$-dimensional Polish space.\smallskip

\noindent (b) Let $t\!\in\!\omega^{<\omega}$. The set $H_t$ is a dense $\bormtwo$ subset of $K^0_t$.\smallskip

\noindent (c) The set $\mathbb{A}_3$ is a $\borathree$ subset of $\mathbb{X}_3^2$. In particular, 
$\mathbb{A}_3$ is an analytic relation on $\mathbb{X}_3$.\smallskip

\noindent (d) Let $\beta\!\in\!\omega^\omega$. Then 
$\bigcap_{n\in\omega}~K_{\beta\vert (n+1)}\!\subseteq\!\mathbb{X}_3$.\smallskip

\noindent (e) $(\mathbb{X}_3,\mathbb{A}_3)\not\preceq_{\borthree}
\big(\omega ,\neg\Delta (\omega )\big)$.\end{lem}

\noindent\bf Proof.\rm\ (a) $\mathbb{X}_3$ is clearly a $\bormtwo$ subset of $2^\omega$. So let us prove its density. We just have to prove that 
$\big\{\alpha\!\in\! 2^\omega\mid\exists p\!\geq\! n~~\alpha\vert p\mbox{ is placed}\big\}$ is dense in 
$2^\omega$ for each integer $n$. So let $\emptyset\!\not=\! w\!\in\! 2^{<\omega}$. Note that 
$\alpha\! :=\! w1^\infty\!\in\! N_{w(0)}\! =\! K^{w(0)}_\emptyset$. Let $p\!\geq\!\mbox{max}(n,\vert w\vert )$ with $(p)_0\! =\! 0$. Then $\alpha\vert (p\! +\! 1)$ is placed with witness $t\! :=\!\emptyset$.\bigskip

\noindent (b) $H_t$ is clearly a $\bormtwo$ subset of $K^0_t$. So let us prove its density. We just have to prove the density in $K^0_t$ of the set $\big\{\alpha\!\in\! K^0_t\mid\exists p\!\geq\! n~~
(p)_0\! =\!\Sigma^t_{\vert t\vert}\mbox{ and }\alpha (p)\! =\! 1\big\}$, for each integer $n$. If 
$\vert t\vert\! =\! 0$, then $K^0_t\! =\! K^0_\emptyset$. As in the proof of (a), we see (with $w(0)\! :=\! 0$) that 
$\alpha\! :=\! w1^\infty\!\in\! N_w\cap K^0_t$ and $p\!\geq\!\mbox{max}(n,\vert w\vert )$ with $(p)_0\! =\! 0$ are suitable. If $\vert t\vert\!\geq\! 1$, then we argue similarly. We put again $w1^\infty$, in the coordinates not determined by $t$.\bigskip

\noindent (c) By (b), $\mathbb{A}_3$ is a $\borathree$ subset of $2^\omega\!\times\! 2^\omega$. So we just have to see that $\mathbb{A}_3\!\subseteq\!\mathbb{X}_3^2$, which is clear.\bigskip

\noindent (d) Let $\alpha\!\in\!\bigcap_{n\in\omega}~K_{\beta\vert (n+1)}$. Note that the sequence 
$\big(\Sigma^{\beta\vert (n+1)}_{n+1}\big)_{n\in\omega}$ is strictly increasing. In particular, 
$p\! :=<\Sigma^{\beta\vert (n+1)}_{n+1},\beta (n\! +\! 1)\! +\! 1>\! +1\!\geq 
<\Sigma^{\beta\vert (n+1)}_{n+1},0>\! +1\!\geq\!\Sigma^{\beta\vert (n+1)}_{n+1}\! +\! 1\!\!\geq\! n$ and 
$t\! :=\! \beta\vert (n+1)$ are witnesses for the fact fact $\alpha\!\in\!\mathbb{X}_3$.\bigskip

\noindent (e) We argue by contradiction, which gives a partition $(C_n)_{n\in\omega}$ of $\mathbb{X}_3$ into $\mathbb{A}_3$-discrete $\bormtwo$ sets. Fix $n\!\in\!\omega$ and $t\!\in\!\omega^{<\omega}$. Let us prove that there is $i\!\in\!\omega$ such that $C_n\cap K_{ti}\! =\!\emptyset$. We argue by contradiction. By Lemma 5.4, $C_n\cap K^\varepsilon_t$ is dense in $K^\varepsilon_t$ for each $\varepsilon\!\in\! 2$. As 
$C_n\cap K^\varepsilon_t$ is $\bormtwo$, it is comeager in $K^\varepsilon_t$. By (b), $H_t$ is also comeager in $K^0_t$, so that this is also the case of $C_n\cap H_t$. In particular, $\varphi_t[C_n\cap H_t]$ is comeager in $K^1_t$, and 
$C_n\cap \varphi_t[C_n\cap H_t]$ too. In particular, $C_n\cap \varphi_t[C_n\cap H_t]$ is not empty, which contradicts the $\mathbb{A}_3$-discreteness of $C_n$.\bigskip

 Applying this inductively, we construct $\beta\!\in\!\omega^\omega$ such that 
$C_n\cap K_{\beta\vert (n+1)}\! =\!\emptyset$ for each $n\!\in\!\omega$. By compactness, there is 
$\alpha\!\in\!\bigcap_{n\in\omega}~K_{\beta\vert (n+1)}$, and 
$\alpha\!\notin\!\bigcup_{n\in\omega}~C_n\! =\!\mathbb{X}_3$. But this contradicts (d).
\hfill{$\square$}\bigskip

 The following uniqueness properties will be important in the sequel.

\begin{lem} Let $t\!\in\! \omega^{<\omega}$ and $\alpha\!\in\! K_t$.\smallskip

\noindent (a) Assume that $u\!\in\! 2^{<\omega}$ is placed with witness $t$. Then 
$\Sigma_t\! <\!\vert u\vert$, the last $1$ in $u$ strictly before the position $<\Sigma^t_{k+1},0>$ is at the position ${<\Sigma^t_k,t(k)\! +\! 1>}$ for each $k\! <\!\vert t\vert$, and $t$ is unique.\smallskip

\noindent (b) Let $p\! >\!\Sigma_t$ be such that $\alpha\vert p$ is placed with witness $t'$. Then 
$t\!\subseteq\! t'$.\end{lem}

\noindent\bf Proof.\rm\ (a) As $(\vert u\vert\! -\! 1)_0\! =\!\Sigma^t_{\vert t\vert}$, we may assume that 
$\vert t\vert\!\geq\! 1$. Let $\alpha\!\in\! N_u\cap K_t$. Then 
$$(\alpha )_{\Sigma^t_k}\! =\!(w_{t(k)})_00^{t(k)+1-\vert (w_{t(k)})_0\vert}10^\infty$$ 
for each $k\! <\!\vert t\vert$.\bigskip

 In particular, $\alpha (<\Sigma^t_k,t(k)\! +\! 1>)\! =\! 1$. As $<n\! +\! p\! +\! 1,0>><n,p>$ and 
$(\vert u\vert\! -\! 1)_0\! =\!\Sigma^t_{\vert t\vert}$, we get $\vert u\vert\! >\!\vert u\vert\! -\! 1\!\geq 
<\Sigma^t_{\vert t\vert},0>=\!\Sigma_t\!\geq <\Sigma^t_{k+1},0>><\Sigma^t_k,t(k)\! +\! 1>$. Thus 
$u(<\Sigma^t_k,t(k)\! +\! 1>)\! =\! 1$.\bigskip

\noindent $\bullet$ Let us prove that the last $1$ in $u$ strictly before the position $<\Sigma^t_{k+1},0>$ is at the position 
${<\Sigma^t_k,t(k)\! +\! 1>}$. The consecutive integers between the values $<\Sigma^t_k,t(k)\! +\! 1>$ and $<\Sigma^t_{k+1},0>$ are ${<\Sigma^t_k,t(k)\! +\! 1>}$, $<\Sigma^t_k\! -\! 1,t(k)\! +\! 2>$, ..., 
$<0,\Sigma^t_k\! +\! t(k)\! +\! 1>$ and $<\Sigma^t_{k+1},0>$. So we have to see that 
$u(<\Sigma^t_k\! -\! j,t(k)\! +\! 1\! +\! j>)\! =\! 0$ if $0\! <\! j\!\leq\!\Sigma^t_k$.

\vfill\eject

 There are $k'\! <\! k$ and $i\! <\! t(k')\! +\! 2$ such that $\Sigma^t_k\! -\! j\! =\!\Sigma^t_{k'}\! +\! i$. In particular, 
$$t(k)\! +\! 1\! +\! j\! =\!\Sigma^t_{k+1}\! -\! 1\! -\!\Sigma^t_{k'}\! -\! i\!\geq\! t(k')\! +\! 2\! -\! i.$$
Lemma 5.2.(c) applied to $w_{t(k')}$ implies that $\vert (w_{t(k')})_i\vert\!\leq\! t(k')\! +\! 2\! -\! i$. Thus 
$$\alpha (<\Sigma^t_{k'}\! +\! i,t(k)\! +\! 1\! +\! j>)\! =\! 0\mbox{,}$$ 
and we are done.\bigskip

\noindent $\bullet$ As the last $1$ in $u$ strictly before the position 
$<\Sigma^t_{\vert t\vert},0>$ is at the position $<\Sigma^t_{\vert t\vert -1},t(\vert t\vert\! -\! 1)\! +\! 1>$, 
$t(\vert t\vert\! -\! 1)$ is determined. It remains to iterate this argument to see the uniqueness of $t$.\bigskip

\noindent (b) We argue by induction on $l\! :=\!\vert t\vert$, and we may assume that our property is proved for $l$. So let $t\!\in\!\omega^{l+1}$, $\alpha\!\in\! K_t$, and $p\! >\!\Sigma_t$ be such that 
$\alpha\vert p$ is placed with witness $t'$. Note that $\alpha\!\in\! K_{t\vert l}$ and 
$p\! >\!\Sigma_{t\vert l}$. By the induction assumption, we get $t\vert l\!\subseteq\! t'$.\bigskip

 Let us prove that $t\vert l\!\not=\! t'$. We argue by contradiction, so that 
$(p\! -\! 1)_0\! =\!\Sigma^{t'}_l\! =\!\Sigma^t_l$. Note that 
${<\Sigma^t_l,t(l)\! +\! 1><\!\Sigma_t\! <\! p}$. Thus $(p\! -\! 1)_1\! >\! t(l)\! +\! 1\! >\! 0$ and 
$\alpha (p\! -\! 1)\! =\! 1\! =\!\alpha (<\Sigma^t_l,(p\! -\! 1)_1>)$. As $\alpha\!\in\! K_t$, the last $1$ of 
$(\alpha )_{\Sigma^t_l}$ is at the position $t(l)\! +\! 1$, which is absurd.\bigskip

This shows the existence of $t'(l)$. Let $\beta\!\in\! N_{\alpha\vert p}\cap K_{t'}$. The last $1$ of 
$(\beta )_{\Sigma^t_l}$ is at the position $t'(l)\! +\! 1$. As 
$<\Sigma^t_l,t'(l)\! +\! 1><\! p\! -\! 1\! =<\Sigma^{t'}_{\vert t'\vert},(p\! -\! 1)_1>$, it is also the last $1$ of 
$(\beta\vert p)_{\Sigma^t_l}\! =\! (\alpha\vert p)_{\Sigma^t_l}$. But the last $1$ of 
$(\alpha\vert p)_{\Sigma^t_l}$ is at the position $t(l)\! +\! 1$, so that $t'(l)\! =\! t(l)$ and $t\!\subseteq\! t'$.
\hfill{$\square$}

\begin{defi} Let $u\!\in\! 2^{<\omega}$ and $l\!\in\!\omega$.\smallskip

\noindent (a) If $u$ is placed, then let $t\!\in\!\omega^{<\omega}$ be the unique witness given by Lemma 5.6.(b). We will consider\smallskip

\noindent $\bullet$ the length $l(u)\! :=\!\vert t\vert$\smallskip

\noindent $\bullet$ the sequence $u^{l(u)}\!\in\! 2^{\vert u\vert}\!\setminus\!\{ u\}$ defined by 
$u^{l(u)}(m)\! :=\! 1\! -\! u(m)$ if $m\! =\!\Sigma_t$, $u(m)$ otherwise. Note that $u^{l(u)}$ is placed with witness $t$, so that $l(u^{l(u)})\! =\! l(u)$ and $(u^{l(u)})^{l(u)}\! =\! u$\smallskip

\noindent $\bullet$ the digit $\varepsilon (u)\! :=\! u(\Sigma_t)$. Note that 
$\varepsilon (u^{l(u)})\! =\! 1\! -\!\varepsilon (u)$.\smallskip

\noindent (b) We say that $u$ is $l$-$placed$ if $u$ is placed and $l(u)\! =\! l$. We say that $u$ is 
$(\leq\! l)$-$placed$ (resp., $(<\! l)$-$placed$, $(>\! l)$-$placed$) if there is $l'\!\leq\! l$ (resp., $l'\! <\! l$, 
$l'\! >\! l$) such that $u$ is $l'$-placed.\end{defi}

 The following lemma will be crucial in the construction of the homomorphism. We construct some finite approximations of the homomorphism. The lemma says that these finite approximations can be constructed independently.

\begin{lem} Let $u\!\not=\! v\!\in\! 2^{<\omega}$ be placed with $\varepsilon (u)\! =\!\varepsilon (v)$. Then 
$\{ u,u^{l(u)}\}\cap\{ v,v^{l(v)}\}\! =\!\emptyset$.\end{lem}

\noindent\bf Proof.\rm\ Note first that 
$\varepsilon (u^{l(u)})\! =\! 1\! -\!\varepsilon (u)\! =\! 1\! -\!\varepsilon (v)\! =\!\varepsilon (v^{l(v)})$. Thus 
$u\!\notin\!\{ v,v^{l(v)}\}$ and $u^{l(u)}\!\not=\! v$. If $u^{l(u)}\! =\! v^{l(v)}$, then 
$u\! =\! (u^{l(u)})^{l(u^{l(u)})}\! =\! (v^{l(v)})^{l(v^{l(v)})}\! =\! v$, which is absurd.\hfill{$\square$}

\vfill\eject

 When we consider the finite approximations of an element of $\mathbb{A}_3$, we have to guess the finite sequence $t$. We usually make some mistakes. In this case, we have to be able to come back to an earlier position. This is the role of the following predecessors.\bigskip

\noindent\bf Notation.\rm\ Let $u\!\in\! 2^{<\omega}$. Note that $<\varepsilon >$ is $0$-placed. This allows us to define
$$u^{-}\! :=\!\left\{\!\!\!\!\!\!\!
\begin{array}{ll}
& \emptyset\mbox{ if }\vert u\vert\!\leq\! 1\mbox{,}\cr
& u\vert\mbox{max}\{ n\! <\!\vert u\vert\mid u\vert n\mbox{ is placed}\}\mbox{ if }\vert u\vert\!\geq\! 2.
\end{array}
\right.$$
and, for $l\!\in\!\omega$,
$$u^{-l}\! :=\!\left\{\!\!\!\!\!\!\!
\begin{array}{ll}
& \emptyset\mbox{ if }\vert u\vert\!\leq\! 1\mbox{,}\cr
& u\vert\mbox{max}\{ n\! <\!\vert u\vert\mid u\vert n\mbox{ is }(\leq\! l)\mbox{-placed}\}\mbox{ if }
\vert u\vert\!\geq\! 2\mbox{,}
\end{array}
\right.$$
Before proving our main theorem, we study the relation between these predecessors and the placed sequences.

\begin{lem} Let $l\!\in\!\omega$ and $u\!\in\! 2^{<\omega}$ be $l$-placed with $\vert u\vert\!\geq\! 2$.\smallskip

\noindent (a) Assume that $u^-$ is $l$-placed. Then $\varepsilon (u^-)\! =\!\varepsilon (u)$. If moreover 
$(u^l)^-$ is $l$-placed, then $(u^l)^-\! =\! (u^-)^l$.\smallskip

\noindent (b) $u^{-l}$ is $l$-placed if and only if $(u^l)^{-l}$ is $l$-placed. In this case, 
$\varepsilon (u^{-l})\! =\!\varepsilon (u)$ and $(u^l)^{-l}\! =\! (u^{-l})^l$.\smallskip

\noindent (c) Assume that $u^-$ or $(u^l)^-$ is $(<\! l)$-placed. Then 
$u^-\! =\! u^{-l}\! =\! (u^l)^-\! =\! (u^l)^{-l}$ is $(l\! -\! 1)$-placed.\smallskip

\noindent (d) Assume that $u^-$ or $(u^l)^-$ is $(>\! l)$-placed. Then exactly one of those two sequences is $(>\! l)$-placed, and the other one is $l$-placed. If $u^-$ (resp., $(u^l)^-$) is $(>\! l)$-placed, then $u^{-l}\! =\!\big( (u^l)^-\big)^l$ (resp., $u^{-l}\! =\! u^-$) and $\varepsilon (u^{-l})\! =\!\varepsilon (u)$ (resp., 
$\varepsilon\big( (u^l)^{-l}\big)\! =\!\varepsilon (u^l)$).\smallskip

\noindent (e) $l(u^{-l})\!\in\!\{ l\! -\! 1,l\}$.\end{lem}

\noindent\bf Proof.\rm\ Let $t\!\in\!\omega^l$ (resp. $t'\!\in\!\omega^{<\omega}$) be a witness for the fact that $u$ (resp., $u^-$) is placed, and $\alpha\!\in\! N_u\cap K_t$.\bigskip

\noindent\bf Claim.\it\ Assume that $(\vert u\vert\! -\! 1)_1\! =\! 0$. Then 
$u^-\! =\! u^{-l}\! =\! (u^l)^-\! =\! (u^l)^{-l}$ is $(l\! -\! 1)$-placed.\bigskip

\noindent\bf Proof.\rm\ Note that $l\!\geq\! 1$ since $\vert u\vert\!\geq\! 2$. The consecutive integers between the values $<\Sigma^t_{l-1},t(l\! -\! 1)\! +\! 1>$ and $\Sigma_t$ are 
${<\Sigma^t_{l-1},t(l\! -\! 1)\! +\! 1>}$, $<\Sigma^t_{l-1}\! -\! 1,t(l\! -\! 1)\! +\! 2>$, ..., 
$<0,\Sigma^t_{l-1}\! +\! t(l\! -\! 1)\! +\! 1>$ and $\Sigma_t$. By Lemma 5.6.(b), $\Sigma_t\! <\!\vert u\vert$ and the last $1$ in $u$ strictly before the position $\Sigma_t$ is at the position 
${<\Sigma^t_{l-1},t(l\! -\! 1)\! +\! 1>}$. This shows that $u\vert (\Sigma_t\! +\! 1)$ and 
$u\vert (<\Sigma^t_{l-1},t(l\! -\! 1)\! +\! 1>+1)$ are placed and 
$\big( u\vert (\Sigma_t\! +\! 1)\big)^-\! =\! u\vert (<\Sigma^t_{l-1},t(l\! -\! 1)\! +\! 1>+1)$ since 
$t(l\! -\! 1)\! +\! j\! >\! 0$ if $1\!\leq\! j\!\leq\!\Sigma^t_{l-1}\! +\! 1$. As $(\vert u\vert\! -\! 1)_1\! =\! 0$, 
$\vert u\vert\! =\!\Sigma_t\! +\! 1$ and the sequence $u^-\! =\! u\vert (<\Sigma^t_{l-1},t(l\! -\! 1)\! +\! 1>\! +1)$ is 
$(l\! -\! 1)$-placed, which implies that $u^-\! =\! u^{-l}\! =\! (u^l)^-\! =\! (u^l)^{-l}$.\hfill{$\diamond$}\bigskip

\noindent (a) By the claim, $(\vert u\vert\! -\! 1)_1\! >\! 0$. Thus $u\vert (\Sigma_t\! +\! 1)\!\subsetneqq\! u$ is $l$-placed, $u\vert (\Sigma_t\! +\! 1)\!\subseteq\! u^-$ and $\Sigma_t\! <\!\vert u^-\vert$. As 
$u^-\!\subseteq\!\alpha$, we can apply Lemma 5.6.(c) and $t\!\subseteq\! t'$. Thus $t\! =\! t'$ since 
$\vert t\vert\! =\!\vert t'\vert\! =\! l$, and the equalities 
$\varepsilon (u^-)\! =\! (u^-)(\Sigma_{t'})\! =\! u(\Sigma_t)\! =\!\varepsilon (u)$ hold.\bigskip

 Assume now that $(u^l)^-$ is $l$-placed. As $u^l$ is $l$-placed with witness $t$, there is some 
$\beta\!\in\! N_{u^l}\cap K_t$. As $u\vert (\Sigma_t\! +\! 1)\!\subseteq\! u^-\!\subsetneqq\! u$, we get 
$\big( u\vert (\Sigma_t\! +\! 1)\big)^l\!\subseteq\! (u^l)^-\!\subset\!\beta$. Thus 
$\vert (u^l)^-\vert\! >\!\Sigma_t$. Lemma 5.6.(c) implies that $t$ is the witness for the fact that $(u^l)^-$ is 
$l$-placed. If $u^-\! =\! u\vert (<\Sigma^t_l,j_0>\! +1)$, then there is no $j_0\! <\! j\! <\! (\vert u\vert\! -\! 1)_1$ with $u(<\Sigma^t_l,j>)\! =\! 1$, and $(u^l)^-\! =\! u^l\vert (<\Sigma^t_l,j_0>\! +1)\! =\! (u^-)^l$.

\vfill\eject

\noindent (b) Assume that $u^{-l}$ is $l$-placed. As in (a) we get $(\vert u\vert\! -\! 1)_1\! >\! 0$ and 
$j_1$ with $u^{-l}\! =\! u\vert (<\Sigma^t_l,j_1>\! +1)$, and 
$(u^l)^{-l}\! =\! u^l\vert (<\Sigma^t_l,j_1>\! +1)\! =\! (u^{-l})^l$ is $l$-placed. The equivalence comes from the fact that $(u^l)^l\! =\! u$. We argue as in (a) to see that $\varepsilon (u^{-l})\! =\!\varepsilon (u)$ if $u^{-l}$ is $l$-placed.\bigskip

\noindent (c) Assume first that $u^-$ is $(<\! l)$-placed. The proof of (a) shows that $\vert t'\vert\!\geq\! l$ if 
$(\vert u\vert\! -\! 1)_1\! >\! 0$. Thus $(\vert u\vert\! -\! 1)_1\! =\! 0$ and the claim gives the result. If $(u^l)^-$  is $(<\! l)$-placed, then we apply this to $u^l$, using the facts that $u^l$ is $l$-placed and $(u^l)^l\! =\! u$.\bigskip

\noindent (d) Assume first that $u^-$ is $(>\! l)$-placed. As in (a) we get $t\!\subsetneqq\! t'$. In particular, the last $1$ in $(u^-)_{\Sigma^t_l}$ is at the position $t'(l)\! +\! 1$. Let us prove that 
$u^{-l}\! =\! u\vert (<\Sigma^t_l,t'(l)\! +\! 1>\! +1)$. Note that $u\vert (<\Sigma^t_l,t'(l)\! +\! 1>\! +1)$ is $l$-placed, so that $u\vert (<\Sigma^t_l,t'(l)\! +\! 1>\! +1)\!\subseteq\! u^{-l}\!\subseteq\! u^-$. Lemma 5.6.(c) shows that $u^{-l}$  is $l$-placed with witness $t$. As the last $1$ in $(u^-)_{\Sigma^t_l}$ is at the position $t'(l)\! +\! 1$, we are done.\bigskip

 Note that $u^l\vert (<\Sigma^t_l,t'(l)\! +\! 1>\! +1)\!\subseteq\! (u^l)^-$. We argue by contradiction to see that $(u^l)^-$ is not $(>\! l)$-placed. This gives a witness $t''$, which is a strict extension of $t$ by Lemma 
5.6.(c). We saw that the last $1$ in $\big( (u^l)^-\big)_{\Sigma^t_l}$ is at the position $t'(l)\! +\! 1$. But it is also at the position $t''(l)\! +\! 1$, which shows that $t''(l)\! =\! t'(l)$. Thus 
$u^l(\Sigma_t)\! =\! w_{t''(l)}(0)\! =\! w_{t'(l)}(0)$. But $u^-(\Sigma_t)\! =\! w_{t'(l)}(0)$. This implies that 
$u^l(\Sigma_t)\! =\! 1\! -\! w_{t'(l)}(0)$, which is absurd. This shows that 
$(u^l)^-\! =\! u^l\vert (<\Sigma^t_l,t'(l)\! +\! 1>\! +1)\! =\! (u^{-l})^l$ is $l$-placed, so that 
$u^{-l}\! =\!\big( (u^l)^-\big)^l$. Moreover, 
$\varepsilon (u^{-l})\! =\! (u^{-l})(\Sigma_t)\! =\! u(\Sigma_t)\! =\!\varepsilon (u)$.\bigskip

 Assume now that $(u^l)^-$ is $(>\! l)$-placed. As $u^l$ is $l$-placed and $(u^l)^l\! =\! l$, the previous arguments show that $u^-$ is $l$-placed. In particular, $u^{-l}\! =\! u^-$.\bigskip

\noindent (e) If $u^-$ is $l$-placed, then $u^{-l}\! =\! u^-$ is $l$-placed. If $u^-$ is $(<\! l)$-placed, then by (c) $u^{-l}$ is $(l\! -\! 1)$-placed. If $u^-$ is $(>\! l)$-placed, then by (d) $(u^l)^-$ is $l$-placed and 
$u^{-l}\! =\!\big( (u^l)^-\big)^l$ is $l$-placed too.\hfill{$\square$}\bigskip

\noindent\bf Proof of Theorem 5.1.\rm\ $\mathbb{X}_3$ and $\mathbb{A}_3$ have been defined before. The ``exactly" part comes from Lemma 5.5.(e). So we just have to prove that (a) or (b) holds. We may assume that $X$ is recursively presented and $A$ is a $\Ana$ relation. We set 
$$U\! :=\!\bigcup\Big\{ V\!\in\!\Ana (X)\cap\bormone\mid\exists D\!\in\!\Borel\cap\bormtwo (\omega\!\times\! X)~~V\!\subseteq\!\bigcup_{p\in\omega}~D_p\mbox{ and }\forall p\!\in\!\omega ~~
A\cap D_p^2\! =\!\emptyset\Big\}.$$
\bf Case 1.\rm\ $U\! =\! X$.\bigskip

 There is a countable covering of $X$ into $A$-discrete $\borathree$ sets. We just have to reduce them to get a partition showing that (a) holds.\bigskip

\noindent\bf Case 2.\rm\ $U\!\not=\! X$.\bigskip

 Note that if $V$ is as in the definition of $U$, then $V$ and $\neg\bigcup_{p\in\omega}~D_p$ are disjoint 
$\Ana$ sets, separable by a $\bormone$ set. By Theorems 1.A and 1.B in [Lo1], there is 
$W\!\in\!\Borel\cap\bormone$ separating these two sets. This shows that we can replace the condition 
``$V\!\in\!\Ana (X)\cap\bormone$" with ``$V\!\in\!\Borel (X)\cap\bormone$" in the definition of $U$. Thus $U$ is $\Ca (X)\cap\boratwo$ since the set of codes for $\Borel\cap\bormxi$ sets is $\Ca$ if 
$\xi\! <\!\omega_1^{\mbox{CK}}$ (see [Lo1]). This shows that $Y\! :=\! X\!\setminus\! U$ is a nonempty 
$\Ana\cap\bormtwo$ subset of $X$.

\vfill\eject

\noindent\bf Claim\it\ If $V\!\in\!\Ana (X)\cap\bormone$ meets $Y$, then $V\cap Y$ is not $A$-discrete.\rm
\bigskip
 
 We argue by contradiction. Note that $V\cap Y\!\in\!\Ana\cap\bormtwo$. Lemma 2.2 gives a 
$\Borel (X)\cap\bormtwo$ set containing $V\cap Y$ and $A$-discrete. As 
$V\cap U$ can be covered with some $\bigcup_{p\in\omega}~D_p$'s, so is $V$. Thus 
$V\!\subseteq\! U$, by $\Borel$-selection. Therefore $V\cap Y\!\subseteq\! U\!\setminus\! U\! =\!\emptyset$, which is absurd.\hfill{$\diamond$}\bigskip

\noindent $\bullet$ We construct, when $u$ is placed, some points $x_u$ of $Y$, some $\Boraone$ subsets $X_u$ of $X$, and some $\Ana$ subsets $U_u$ of $X^2$. We want these objects to satisfy the following conditions:
$$\begin{array}{ll}
& (1)~x_u\!\in\! X_u\mbox{ and}\left\{\!\!\!\!\!\!\!\!
\begin{array}{ll}
& (x_u,x_{u^{l(u)}})\!\in\! U_u\mbox{ if }\varepsilon (u)\! =\! 0\cr\cr
& (x_{u^{l(u)}},x_u)\!\in\! U_u\mbox{ if }\varepsilon (u)\! =\! 1
\end{array}
\right.\cr
& (2)~\overline{X_u}\!\subseteq\! X_{u^-}\mbox{ if }\vert u\vert\!\geq\! 2\cr
& (3)~U_u\! =\! U_{u^{l(u)}}\!\subseteq\! A\cap Y^2\cap\Omega_{X^2}\mbox{, and }
U_u\!\subseteq\! U_{u^{-l}}\mbox{ if }u\mbox{ and }u^{-l}\mbox{ are }l\mbox{-placed}\cr
& (4)~\mbox{diam}(X_u)\!\leq\! 2^{-\vert u\vert}\mbox{ and }
\mbox{diam}_{\mbox{GH}}(U_u)\!\leq\! 2^{-\vert u\vert}\cr
& (5)~U_u\!\!\subseteq\!\!\left\{\!\!\!\!\!\!\!\!
\begin{array}{ll}
& \overline{\Pi_0[(X_{u^{-l}}\!\times\! X_{(u^{-l})^{l(u^{-l})}})\cap U_{u^{-l}}]}^2
\mbox{ if }\varepsilon (u^{-l})\! =\! 0\cr\cr
& \overline{\Pi_1[(X_{(u^{-l})^{l(u^{-l})}}\!\times\! X_{u^{-l}})\cap U_{u^{-l}}]}^2
\mbox{ if }\varepsilon (u^{-l})\! =\! 1
\end{array}
\right.
\!\!\mbox{if }u\mbox{ is }l\mbox{-placed and }u^{-l}\mbox{ is not}
\end{array}$$

 As we will see, Conditions (1)-(4) are sufficient to get the required objects. Condition (5) is used to prove that the construction is possible. The idea is the following. When we extend some $u\!\in\! 2^{<\omega}$, some new links may appear. But we may also break some links, and preserve only an initial segment of them. In this case, to ensure Condition (3), we have to be able to come back to the last preserved link. This is possible if we use iteratively Conditions (3) and (5).\bigskip

\noindent $\bullet$ Assume that this is done.  Let $\alpha\!\in\!\mathbb{X}_3$ and 
$(p^\alpha_k)_{k\in\omega}$ be the infinite strictly increasing sequence of integers $p^\alpha_k\!\geq\! 1$ such that $\alpha\vert p^\alpha_k$ is placed. Note that 
$\overline{X_{\alpha\vert p^\alpha_{k+1}}}\!\subseteq\! X_{\alpha\vert p^\alpha_k}$, by Condition (2). This shows that $(\overline{X_{\alpha\vert p^\alpha_k}})_{k\in\omega}$ is a non-increasing sequence of nonempty closed subsets of $X$ whose diameters tend to $0$, and we define 
$\{ f(\alpha )\}\! :=\!\bigcap_{k\in\omega}~\overline{X_{\alpha\vert p^\alpha_k}}\! =\!
\bigcap_{k\in\omega}~X_{\alpha\vert p^\alpha_k}$, so that $f\! :\!\mathbb{X}_3\!\rightarrow\! X$ is continuous and $f(\alpha )\! =\!\mbox{lim}_{k\rightarrow\infty}~x_{\alpha\vert p^\alpha_k}$.\bigskip

 Now let $(\alpha ,\beta )\!\in\!\mathbb{A}_3$. If 
$(\alpha ,\beta )\!\in\!\bigcup_{\vert t\vert =l}~\mbox{Gr}({\varphi_t}_{\vert H_t})$, then let 
$(p_j)_{j\in\omega}$ be the infinite strictly increasing sequence of integers $p_j\!\geq\! 1$ such that 
$(p_j\! -\! 1)_0\! =\!\Sigma^t_{\vert t\vert}$, $(p_j\! -\! 1)_1\! >\! 0$ and $\alpha (p_j\! -\! 1)\! =\! 1$. In particular, 
$\alpha\vert p_j$ is $l$-placed and $\varepsilon (\alpha\vert p_j)\! =\! 0$. Note that $(p_j)_{j\in\omega}$ is also the infinite strictly increasing sequence of integers $p_j\!\geq\! 1$ such that 
$(p_j\! -\! 1)_0\! =\!\Sigma^t_{\vert t\vert}$, $(p_j\! -\! 1)_1\! >\! 0$ and $\beta (p_j\! -\! 1)\! =\! 1$ on one side, and a subsequence of both $(p^\alpha_k)_{k\in\omega}$ and $(p^\beta_k)_{k\in\omega}$ on the other side.\bigskip

 If moreover $p\!\geq\! p_0$ and $\alpha\vert p$ is placed, then the witness is an extension of $t$ and 
$l(\alpha\vert p)\!\geq\! l$, by Lemma 5.6.(c). In particular, if $p\!\geq\! p_0$ and $\alpha\vert p$ is $l$-placed, then the witness is $t$. This proves that $(p_j)_{j\in\omega}$ is the infinite strictly increasing sequence of integers $p_j\!\geq\! p_0$ such that $\alpha\vert p_j$ is $l$-placed. Therefore 
$(\alpha\vert p_{j+1})^{-l}\! =\!\alpha\vert p_j$. By Condition (3), 
$(U_{\alpha\vert p_j})_{j\in\omega}$ is a non-increasing sequence of nonempty clopen subsets of 
$A\cap\Omega_{X^2}$ whose GH-diameter tend to $0$. So we can define $F(\alpha ,\beta )\!\in\! A$ by 
$\{ F(\alpha ,\beta )\}\! :=\!\bigcap_{j\in\omega}~U_{\alpha\vert p_j}$. Note that 
$F(\alpha ,\beta )\! =\!\mbox{lim}_{j\rightarrow\infty}~(x_{\alpha\vert p_j},x_{\beta\vert p_j})\! =\!
\big( f(\alpha ),f(\beta )\big)\!\in\! A$, so that $\mathbb{A}_3\!\subseteq\! (f\!\times\! f)^{-1}(A)$.\bigskip

\noindent $\bullet$ Let us prove that the construction is possible. We do it by induction on the length $k$ of $u$.\bigskip

\noindent\bf Subcase 1. $k\! =\! 0$\rm\bigskip

 We are done since $\emptyset$ is not placed.\bigskip

\noindent\bf Subcase 2. $k\! =\! 1$\rm\bigskip

 The claim gives $(x_0,x_1)\!\in\! A\cap Y^2\cap\Omega_{X^2}$. We choose a $\Boraone$ neighborhood $X_{\varepsilon'}$ of $x_{\varepsilon'}$ with diameter at most $2^{-1}$, as well as a 
$\Ana$ subset $U_0\! =\! U_1$ of $X^2$ with GH-diameter at most $2^{-1}$ such that 
$(x_0,x_1)$ is in $U_0\!\subseteq\! A\cap Y^2\cap\Omega_{X^2}$. We are done since 
$<\varepsilon'>$ is $0$-placed and $\varepsilon (<\varepsilon'>)\! =\!\varepsilon'$.\bigskip

\noindent\bf Subcase 3. $k\!\geq\! 2$\rm\bigskip

 If there is no placed sequence in $2^k$, then there is nothing to do. If $u\!\in\! 2^k$ is $l$-placed, then 
$u^l\!\in\! 2^k$ is $l$-placed and $\varepsilon (u^l)\! =\! 1\! -\!\varepsilon (u)$. Assume for example that 
$\varepsilon (u)\! =\! 0$. Lemma 5.8 ensures that we just have to define $x_u,x_{u^l},X_u,X_{u^l}$ and $U_u\! =\! U_{u^l}$, independently from the other sequences in $2^k$.\bigskip

\noindent - If $u^-$ and $(u^l)^-$ are $l$-placed, then $\big( u^-,(u^l)^-\big)\! =\!\big( u^{-l},(u^l)^{-l}\big)$. Moreover, $\varepsilon (u^-)\! =\!\varepsilon (u)\! =\! 0$, $(u^l)^-\! =\! (u^-)^l$ and $(u^l)^{-l}\! =\! (u^{-l})^l$, by Lemma 5.9. We set $(x_u,x_{u^l})\! :=\! (x_{u^-},x_{(u^-)^l})$, we choose $\Boraone$ sets $X_u,X_{u^l}$ with diameter at most $2^{-k}$ such that $(x_u,x_{u^l})\!\in\! X_u\!\times\! X_{u^l}\!\subseteq\!
\overline{X_u}\!\times\!\overline{X_{u^l}}\!\subseteq\! X_{u^-}\!\times\! X_{(u^-)^l}$, as well as a $\Ana$ subset $U_u$ of $X^2$ with GH-diameter at most $2^{-k}$ such that 
$(x_u,x_{u^l})\!\in\! U_u\!\subseteq\! U_{u^-}\! =\! U_{u^{-l}}$. We are done since 
$U_{u^l}\! =\! U_u\!\subseteq\! U_{u^-}\! =\! U_{(u^-)^l}\! =\! U_{(u^{-l})^l}\! =\! U_{(u^l)^{-l}}$.\bigskip

\noindent - If $u^-$ or $(u^l)^-$ is $(<\! l)$-placed, then $u^-\! =\! u^{-l}\! =\! (u^l)^-\! =\! (u^l)^{-l}$ is 
$(l\! -\! 1)$-placed, by Lemma 5.9.(c). Let $W$ be a $\Boraone$ neighborhood of 
$x_{u^-}\! =\! x_{(u^l)^-}$ with $\overline{W}\!\subseteq\! X_{u^-}$. Note that  
$$x_{u^-}\!\in\!\left\{\!\!\!\!\!\!\!
\begin{array}{ll}
& \overline{\Pi_0[(X_{u^-}\!\times\! X_{(u^-)^{l-1}})\cap U_{u^-}]}
\mbox{ if }\varepsilon (u^-)\! =\! 0\mbox{,}\cr\cr
& \overline{\Pi_1[(X_{(u^-)^{l-1}}\!\times\! X_{u^-})\cap U_{u^-}]}
\mbox{ if }\varepsilon (u^-)\! =\! 1.
\end{array}
\right.$$
Assume for example that we are in the second case. Then 
$$x_{u^-}\!\in\! 
\overline{W}\cap\overline{\Pi_1[(X_{(u^-)^{l-1}}\!\times\! X_{u^-})\cap U_{u^-}]}\cap Y
\!\not=\!\emptyset .$$

 The claim gives a couple $(x_u,x_{u^l})\!\in\! A\cap 
(\overline{W}\cap\overline{\Pi_1[(X_{(u^-)^{l-1}}\!\times\! X_{u^-})\cap U_{u^-}]}\cap Y)^2\cap
\Omega_{X^2}$ since the set 
$\overline{W}\cap\overline{\Pi_1[(X_{(u^-)^{l-1}}\!\times\! X_{u^-})\cap U_{u^-}]}$ is 
$\Ana\cap\bormone$. We choose $\Boraone$ sets $X_u,X_{u^l}$ with diameter at most $2^{-k}$ such that 
$(x_u,x_{u^l})\!\in\! X_u\!\times\! X_{u^l}\!\subseteq\!\overline{X_u}\!\times\!\overline{X_{u^l}}\!\subseteq\! 
X_{u^-}\!\times\! X_{(u^l)^-}$, as well as a $\Ana$ subset $U_u$ of $X^2$ with GH-diameter at most 
$2^{-k}$ such that $(x_u,x_{u^l})\!\in\! U_u\!\subseteq\! 
A\cap\overline{\Pi_1[(X_{(u^-)^{l-1}}\!\times\! X_{u^-})\cap U_{u^-}]}^2\cap Y^2\cap\Omega_{X^2}$.\bigskip

\noindent - If $u^-$ or $(u^l)^-$ is $(>\! l)$-placed, then by Lemma 5.9.(d) exactly one of those two sequences is $(>\! l)$-placed, and the other one is $l$-placed. If $u^-$ (resp., $(u^l)^-$) is $(>\! l)$-placed, then $u^{-l}\! =\!\big( (u^l)^-\big)^l$ (resp., $u^{-l}\! =\! u^-$). So assume first that 
$u^-$ is $(>\! l)$-placed, so that $(u^l)^-\! =\! (u^l)^{-l}\! =\! (u^{-l})^l$ and $u^{-l}$ is $l$-placed. Here is an illustration of what is going on in this case.\bigskip

\includegraphics[scale=0.713]{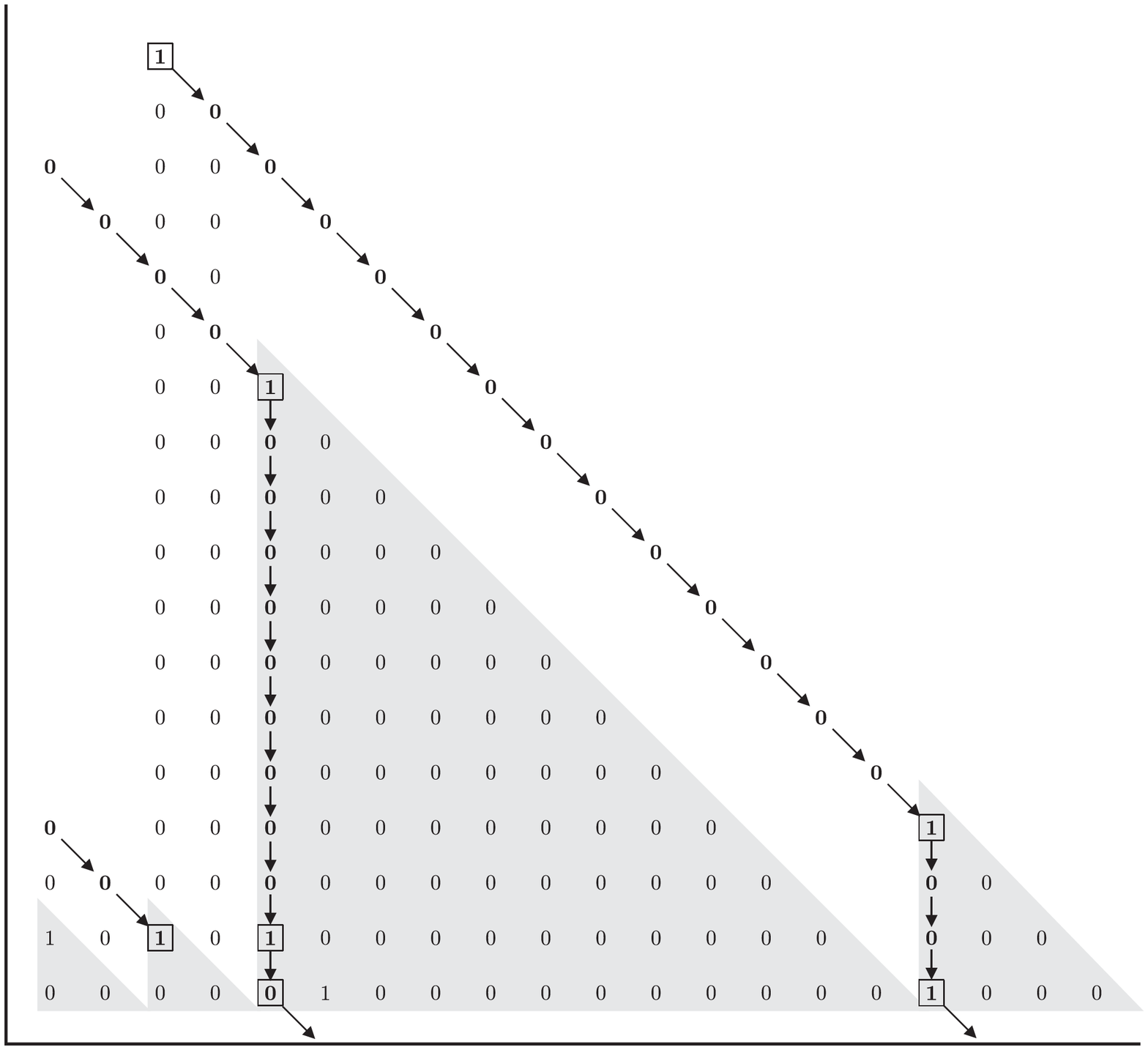}\bigskip

 We define $(u_n)_{n\leq L}$ by $u_0\! :=\! u^-$, $u_L\! :=\! u^{-l}$ and $u_{n+1}\! :=\! u_n^{-l(u_n)}$ if 
$n\! <\! L$. This can be done, by Lemma 5.9.(e). Note that $u_n$ is placed. We enumerate injectively the sequence $\big( l(u_n)\big)_{{n\leq L}}$ by the non-increasing sequence 
$(l_k)_{k\leq K}$. More concretely, ${K\! =\! l_0\! -\! l\!\geq\! 1}$, 
$l(u_0)\! =\! ...\! =\! l(u_{N_0-1})\! =\! l_0$, $l(u_{N_0})\! =\! ...\! =\! l(u_{N_0+N_1-1})\! =\! l_1\! =\! l_0\! -\! 1$, ..., 
$$l(u_{N_0+...+N_{K-2}})\! =\! ...\! =\! l(u_{N_0+...+N_{K-1}-1})\! =\! l_{K-1}\! =\! l_{K-2}\! -\! 1\! =\! l\! +\! 1$$ 
and $l(u_{N_0+...+N_{K-1}})\! =\! l(u_L)\! =\! l_K\! =\! l$, with $N_0,...,N_{K-1}\!\geq\! 1$.\bigskip 

 Note that $u_1\! =\! u_0^{-l_0}$ is $l_0$-placed if $N_0\!\geq\! 2$. By Condition (3), we get 
$U_{u_0}\!\subseteq\! U_{u_1}$. We can iterate this argument, so that the inclusion 
$U_{u_0}\!\subseteq\! U_{u_{N_0-1}}$ holds, even if $N_0\! =\! 1$. By Condition (5), $x_{u^-}$ is in 
$$X_{u^-}\cap\left\{\!\!\!\!\!\!
\begin{array}{ll}
& \overline{\Pi_0[(X_{u_{N_0}}\!\times\! X_{(u_{N_0})^{l_1}})\cap U_{u_{N_0}}]}\mbox{ if }
\varepsilon (u_{N_0})\! =\! 0\mbox{,}\cr\cr
& \overline{\Pi_1[(X_{(u_{N_0})^{l_1}}\!\times\! X_{u_{N_0}})\cap U_{u_{N_0}}]}\mbox{ if }
\varepsilon (u_{N_0})\! =\! 1.
\end{array}
\right.$$ 
Similarly, $X_{u^-}\!\subseteq\! X_{u_{N_0}}$.

\vfill\eject

 This gives 
$$(x^0_0,x^0_1)\!\in\!\left\{\!\!\!\!\!\!
\begin{array}{ll}
& (X_{u^-}\!\times\! X_{(u_{N_0})^{l_1}})\cap U_{u_{N_0}}\mbox{ if }\varepsilon (u_{N_0})\! =\! 0\mbox{,}\cr\cr
& (X_{(u_{N_0})^{l_1}}\!\times\! X_{u^-})\cap U_{u_{N_0}}\mbox{ if }\varepsilon (u_{N_0})\! =\! 1.
\end{array}
\right.$$
If we iterate the previous argument, then we get
$$(x_u,x_{u^l})\! :=\! (x^{K-1}_0,x^{K-1}_1)\!\in\! (X_{u^-}\!\times\! X_{(u_{N_0+...+N_{K-1}})^{l_K}})\cap 
U_{u_{N_0+...+N_{K-1}}}\! =\! (X_{u^-}\!\times\! X_{(u^{-l})^l})\cap U_{u^{-l}}$$
since $\varepsilon (u^{-l})\! =\!\varepsilon (u)\! =\! 0$. We choose $\Boraone$ sets $X_u,X_{u^l}$ with diameter at most $2^{-k}$ such that 
$$(x_u,x_{u^l})\!\in\! X_u\!\times\! X_{u^l}\!\subseteq\!\overline{X_u}\!\times\!\overline{X_{u^l}}\!\subseteq\! 
X_{u^-}\!\times\! X_{(u^{-l})^l}\! =\! X_{u^-}\!\times\! X_{(u^l)^-}\mbox{,}$$ 
as well as a $\Ana$ subset $U_u$ of $X^2$ with GH-diameter at most $2^{-k}$ such that 
$(x_u,x_{u^l})\!\in\! U_u\!\subseteq\! U_{u^{-l}}$.\bigskip

 If now $(u^l)^-$ is $(>\! l)$-placed, then we argue similarly, using the fact that 
$$\varepsilon\big( (u^l)^{-l}\big)\! =\!\varepsilon (u^l)\! =\! 1\! -\!\varepsilon (u)\! =\! 1.$$
This finishes the proof.\hfill{$\square$}\bigskip

 At the beginning of the section, we mentioned the fact that it is not necessary to use the dense $G_\delta$ subset $\mathbb{X}_3$ of $2^\omega$ to find a relation $\mathbb{G}$ on $2^\omega$ satisfying 
$(2^\omega ,\mathbb{G})\not\preceq_{\borthree}\big(\omega ,\neg\Delta (\omega )\big)$. We now specify this.\bigskip

\noindent\bf Notation.\rm\ We set, for $t\!\in\!\omega^{<\omega}$, $\tilde H_t\! :=\! 
K^0_t\!\setminus\!\big(\bigcup_{n\in\omega ,w_n(0)=0}~K_{tn}\cup
\bigcup_{n\in\omega ,w_n(0)=1}~\varphi_t^{-1}(K_{tn})\big)$. Note that $\tilde H_t$ is a $\bormtwo$ subset of $2^\omega$ and $\tilde H_t\cap\varphi_t[\tilde H_t]\! =\!\emptyset$.

\begin{lem} (a) The $\tilde H_t\cup\varphi_t[\tilde H_t]$'s are pairwise disjoint.\smallskip

\noindent (b) The set $H_t$ is a subset of $\tilde H_t$, and thus satisfies the previous disjointness properties.\end{lem}

\noindent\bf Proof.\rm\ (a) Note first that $K_{tn}\!\subseteq\! K^{w_n(0)}_t\!\subseteq\! K_t$ and 
$(K_{tn})_n$ is a sequence of pairwise disjoint sets. This implies that $K_t\cap K_{t'}\! =\!\emptyset$ if $t,t'$ are incompatible. In particular, as 
$\tilde H_t\!\subseteq\! K^0_t\!\setminus\! (\bigcup_{n\in\omega ,w_n(0)=0}~K_{tn})\!\subseteq\! K_t$ and 
$\varphi_t[\tilde H_t]\!\subseteq\! K^1_t\!\setminus\! (\bigcup_{n\in\omega ,w_n(0)=1}~K_{tn})\!\subseteq\! 
K_t$, we also get 
$(\tilde H_t\cup\varphi_t[\tilde H_t])\cap (\tilde H_{t'}\cup\varphi_{t'}[\tilde H_{t'}])\! =\!\emptyset$ if $t$ and $t'$ are incompatible. Now 
$$\begin{array}{ll}
\tilde H_t\cap\tilde H_{tt'n}\!\!
& \!\subseteq\!\tilde H_t\cap K_{t[(t'n)\vert 1]}\cr\cr
& \!\subseteq\!\left\{\!\!\!\!\!\!\!\!
\begin{array}{ll}
& K^0_t\cap K^1_t\mbox{ if }s_{(t'n)(0)}(0)\! =\! 1\mbox{,}\cr
& \neg K_{t[(t'n)\vert 1]}\cap K_{t[(t'n)\vert 1]}\mbox{ if }s_{(t'n)(0)}(0)\! =\! 0
\mbox{,}
\end{array}
\right.
\end{array}$$
so that $\tilde H_t\cap\tilde H_{t'}\! =\!\emptyset$ if $t\!\not=\! t'$. Similarly, 
$(\tilde H_t\cup\varphi_t[\tilde H_t])\cap (\tilde H_{t'}\cup\varphi_{t'}[\tilde H_{t'}])\! =\!\emptyset$ if 
$t\!\not=\! t'$.\bigskip

\noindent (b) If $\alpha\!\in\! K_{tn}\cup\varphi_t^{-1}(K_{tn})$, then $(\alpha )_{\Sigma^t_{\vert t\vert}}$ has finitely many $1$'s.\hfill{$\square$}\bigskip

\noindent\bf Remarks.\rm\ (a) We set 
$\mathbb{G}\! :=\!\bigcup_{t\in\omega^{<\omega}}~\mbox{Gr}({\varphi_t}_{\vert\tilde H_t})$, so that 
$(2^\omega ,\mathbb{G})\not\preceq_{\borthree}\big(\omega ,\neg\Delta (\omega )\big)$, by the proof of Lemma 5.5. By Lemma 5.10, $\mathbb{G}$ is the Borel graph of a partial injection, as announced at the beginning of the section.\bigskip

\noindent (b) Note that 
$(\mathbb{X}_3,\mathbb{A}_3)\preceq_{\borfour}\big(\omega ,\neg\Delta (\omega )\big)$, as we can see with the following partition of $\mathbb{X}_3$:
$$\mathbb{X}_3\! =\!\bigcup_{t\in\omega^{<\omega}}~H_t\cup
\big( \mathbb{X}_3\!\setminus\! (\bigcup_{t\in\omega^{<\omega}}~H_t)\big)\mbox{,}$$
with $H_t\!\in\!\bormtwo$ and $\mathbb{A}_3$-discrete by Lemma 5.10.(b), 
$\mathbb{X}_3\!\setminus\! (\bigcup_{t\in\omega^{<\omega}}~H_t)\!\in\!\bormthree$ and $\mathbb{A}_3$-discrete.\bigskip

\noindent (c) There are a comparing $2$-disjoint family 
$(C^\varepsilon_i)_{(\varepsilon ,i)\in 2\times\omega}$ of subsets of $\mathbb{X}_3$, and also homeomorphisms $\varphi_i\! :\! C^0_i\!\rightarrow\! C^1_i$ such that 
$\mathbb{A}_3\! =\!\bigcup_{i\in\omega}~\mbox{Gr}(\varphi_i)$. Indeed, we choose a bijection 
$b\! :\!\omega\!\rightarrow\!\omega^{<\omega}$ with $b^{-1}(s)\!\leq\! b^{-1}(t)$ if 
$s\!\subseteq\! t$, as in the proof of Lemma 4.3, and set $C^0_i\! :=\! H_{b(i)}$, 
$C^1_i\! :=\!\varphi_{b(i)}[H_{b(i)}]$, $\varphi_i\! :=\! {\varphi_{b(i)}}_{\vert H_{b(i)}}$, so that 
$\mathbb{A}_3\! =\!\bigcup_{i\in\omega}~\mbox{Gr}(\varphi_i)$. It remains to see that 
$(C^\varepsilon_i)_{(\varepsilon ,i)\in 2\times\omega}$ is comparing. We set
$$O^p_q\! :=\!\left\{\!\!\!\!\!\!\!
\begin{array}{ll}
& K_{b(i)}^\varepsilon\!\setminus\! (\bigcup_{l\leq q,b(i)n\subseteq b(l), w_n(0)=\varepsilon}~K_{b(i)n})
\mbox{ if }p\! =\! 2i\! +\!\varepsilon\!\leq\! 2q\! +\! 1\mbox{,}\cr
& \mathbb{X}_3\!\setminus\! (\bigcup_{p'\leq 2q+1}~O^{p'}_q)\mbox{ if }p\! =\! 2q\! +\! 2\mbox{,}\cr
& \emptyset\mbox{ if }p\!\geq\! 2q\! +\! 3\mbox{,}
\end{array}
\right.$$
so that $(O_q^p)_{p\in\omega}$ is a partition of $\mathbb{X}_3$ into $\bortwo$ sets since 
$K_{tn}\!\subseteq\! K^{w_n(0)}_t\!\subseteq\! K_t$, 
$K^0_t\cap K^1_t\! =\!\emptyset$ and $K_t\cap K_{t'}\! =\!\emptyset$ if $t$ and $t'$ are incompatible.\bigskip

 As $\tilde H_t\!\subseteq\! K^0_t\!\setminus\! (\bigcup_{n\in\omega ,w_n(0)=0}~K_{tn})$,  
$\varphi_t[\tilde H_t]\!\subseteq\! K^1_t\!\setminus\! (\bigcup_{n\in\omega ,w_n(0)=1}~K_{tn})$ and 
$H_t\!\subseteq\!\tilde H_t$, (b) in Definition 4.2 is fulfilled. If $q\! <\! i$, then\bigskip

\noindent - either there is no $j\!\leq q$ such that $b(i)$ is compatible with $b(j)$. 
$C^0_i\cup C^1_i\!\subseteq\! K_{b(i)}\!\subseteq\! O_q^{2q+2}$ and we set 
$p_q^i\! :=\! 2q\! +\! 2$.\bigskip

\noindent - or there are $j\!\leq q$ and $n$ such that 
$b(j)n\!\subseteq\! b(i)$, in which case $K^{}_{b(i)}\!\subseteq\! K^{w_n(0)}_{b(j)}\cap K^{}_{b(j)n}$. In particular, $K_{b(i)}$ is disjoint from or included in each difference 
$K^{w_{n'}(0)}_{b(j')}\!\setminus\! K_{b(j')n'}$. Thus $K_{b(i)}$ is disjoint from or included in 
$O^{2j+\varepsilon}_q$. By disjointness, there is at most one couple 
$(j,\varepsilon )$ such that $K_{b(i)}\!\subseteq\! O^{2j+\varepsilon}_q$. If it exists, then we set 
$p_q^i\! :=\! 2j\! +\!\varepsilon$. If it does not exist, then we set $p_q^i\! :=\! 2q\! +\! 2$.

\section{$\!\!\!\!\!\!$ References}

\noindent [K]\ \ A. S. Kechris,~\it Classical Descriptive Set Theory,~\rm 
Springer-Verlag, 1995

\noindent [K-S-T]\ \ A. S. Kechris, S. Solecki and S. Todor\v cevi\'c, Borel chromatic numbers,\ \it 
Adv. Math.\rm\ 141 (1999), 1-44

\noindent [L1]\ \ D. Lecomte, On minimal non potentially closed subsets of the plane,
\ \it Topology Appl.~\rm 154, 1 (2007) 241-262

\noindent [L2]\ \ D. Lecomte, A dichotomy characterizing analytic graphs of uncountable Borel chromatic number in any dimension,~\it Trans. Amer. Math. Soc.\rm~361 (2009), 4181-4193

\noindent [L3]\ \ D. Lecomte, How can we recognize potentially $\bormxi$ subsets of the plane?,
~\it J. Math. Log.\rm~9, 1 (2009), 39-62

\noindent [Lo1]\ \ A. Louveau, A separation theorem for $\Ana$ sets,\ \it Trans. 
Amer. Math. Soc.\ \rm 260 (1980), 363-378

\noindent [Lo2]\ \ A. Louveau, Ensembles analytiques et bor\'eliens dans les 
espaces produit,~\it Ast\'erisque (S. M. F.)\ \rm 78 (1980)

\noindent [M]\ \ Y. N. Moschovakis,~\it Descriptive set theory,~\rm North-Holland, 1980
 
\end{document}